\documentclass[letterpaper,12pt]{article}

\usepackage{amsmath}
\allowdisplaybreaks
\usepackage{float}
\usepackage{microtype,booktabs}
\usepackage{amsmath,amssymb,amsfonts,amsthm,bm}
\usepackage{authblk}
\usepackage{mathrsfs,dsfont}
\usepackage{indentfirst}
\usepackage[pdftex]{color,graphicx}
\usepackage[pdftex,bookmarks,unicode,colorlinks,linkcolor=blue]{hyperref}
\usepackage{tikz}
\usepackage{makecell}

\usepackage{graphicx}
\usepackage{subcaption}
\usetikzlibrary{arrows, decorations.pathmorphing, backgrounds, positioning, fit, petri, automata,cd}
\usepackage{rotating}
\numberwithin{equation}{section}

\newtheorem{theorem}{Theorem}[section] 
     
\newtheorem{lemma}{Lemma}[section]       
\newtheorem{pro}{Proposition}[section]   
\theoremstyle{remark}
\newtheorem{remark}{Remark}[section]     
\theoremstyle{definition}
\newtheorem{definition}{Definition}[section] 
   
\theoremstyle{definition}
\newtheorem{assumption}{Assumption}[section]

\title{Onsager--Machlup Theory and Transition Pathways for Degenerate Stochastic Systems with Multiplicative Noise}

\date{}
\author[$\dag$]{\small Weichen Wu}
\author[$\dag$]{\small Ying Chao\thanks{Corresponding author\par \textit{Email:} yingchao1993@xjtu.edu.cn}}
\author[$\ddag$]{\small Jinqiao Duan}

\affil[$\dag$]{\footnotesize School of Mathematics and Statistics, Xi'an Jiaotong University, Xi'an, Shaanxi 710049, China.}
\affil[$\ddag$]{Department of Mathematics \& Guangdong Provincial Key Laboratory of Mathematical and Neural Dynamical Systems, Great Bay University, Dongguan, Guangdong 523000, China}

\begin{document}

\maketitle

\begin{abstract}
Noise-induced transitions between metastable states are central to many
nonlinear stochastic dynamical systems, and most probable transition pathways
provide an effective means of revealing the underlying transition mechanisms.
In this work, we develop an Onsager--Machlup framework for degenerate stochastic
dynamical systems driven by multiplicative Gaussian noise. We first derive the
Onsager--Machlup action functional for the additive-noise case using the
Girsanov transformation, and then extend the result to multiplicative noise
through a Lamperti-type transformation. The resulting functional provides a
variational characterization of the most probable transition pathway, from which we derive the associated
Hamilton--Pontryagin boundary value problem. We further recast the variational
problem as an equivalent optimal control problem and propose a
neural-network-based numerical method for computing the most probable transition
pathways. Numerical examples are presented to illustrate the
theoretical framework and demonstrate the effectiveness of the proposed
computational approach.
\medskip


\end{abstract}

\section{Introduction}
Stochastic dynamical systems provide a fundamental framework for modeling
complex systems subject to random perturbations. In contrast to deterministic
dynamics, stochastic systems may exhibit noise-induced behaviors that cannot be
inferred from the deterministic flow alone. One important manifestation is
metastability, where trajectories remain near a long-lived state for a
prolonged period before random fluctuations trigger a transition to another
metastable state. Such noise-induced transitions arise in a wide range of
natural systems, including gene regulatory systems, where random fluctuations may induce switches between high and low concentrations of
transcriptional activators
\cite{cheng_most_2019,zheng_transitions_2016,liu_fluctuations-induced_2004},
and carbon-cycle dynamics, where transitions can be reflected by changes in
inorganic carbon and carbonate-ion concentrations \cite{wei_most_2021}.
Understanding the mechanisms of these transitions and identifying their most
probable pathways remain significant challenges.

The most probable transition pathway (MPTP) is a central concept in the study
of noise-induced transitions. It describes a path along which a stochastic
system is most likely to move between metastable states. Formally, the
likelihood of sample paths remaining in a small neighborhood of a prescribed
smooth curve can be asymptotically characterized by an action functional on
path space. The Onsager--Machlup (OM) functional provides such an action
functional and serves as a criterion for comparing the relative likelihoods of
admissible paths. This motivates the use of the OM functional to characterize
and compute MPTPs.

The OM functional was introduced by Onsager and Machlup in their study of
thermodynamic fluctuations \cite{onsager_fluctuations_1953}, originally for
systems with linear drift and constant diffusion coefficients. Tisza and
Manning later provided an axiomatic formulation \cite{tisza_fluctuations_1957},
and Stratonovich developed a probabilistic approach to diffusion path
probabilities \cite{stratonovich1971probability}. Subsequent works refined the
derivation of OM functionals under various assumptions
\cite{durr_onsager-machlup_1978,horsthemke_onsager-machlup_1975,bach_functionals_1977,capitaine_onsager-machlup_1995},
with extensions also made to systems driven by L\'evy noise
\cite{chao_onsagermachlup_2019,hu_transition_2021,WOS:001473901800007,WOS:001514271800001}.
For Brownian diffusions, however, most existing OM theories assume a
nondegenerate diffusion coefficient and therefore do not apply directly to
degenerate stochastic dynamical systems. Although degenerate additive noise
with a constant diffusion coefficient has been studied
\cite{aihara_mortensen_1999,liu_onsagermachlup_2024}, degenerate stochastic
dynamical systems with multiplicative Gaussian noise remain substantially more
challenging, since the state dependence and degeneracy of the diffusion
coefficient must be handled simultaneously. This setting is
natural in applications \cite{baxendale_lyapunov_2003,zhan_numerical_2023} where random fluctuations act only on selected
components of the system and where the noise intensity depends on the current
state. In such cases, transition pathways are shaped not only by the drift
field but also by the geometry induced by the degenerate, state-dependent
diffusion coefficient. This gap motivates the present study.

We also investigate the characterization and computation of MPTPs for
degenerate stochastic systems with multiplicative Gaussian noise. Within the OM
framework, two computational strategies are particularly relevant to the
present work. The first treats the OM functional as a Lagrangian, derives the
associated Euler--Lagrange equations from a variational principle, and solves
the resulting boundary value problem by a shooting method
\cite{cheng_most_2019,chao_onsagermachlup_2019}. However, shooting-type
methods generally do not scale well to high-dimensional problems
\cite{li2021machine,soons_optimal_2024}. Moreover, in degenerate systems, the degeneracy of
the diffusion coefficient introduces constraints into the variational
formulation and leads to more complicated Hamilton--Pontryagin equations
\cite{WOS:001514271800001}. The second strategy formulates the computation of
MPTPs as a constrained optimal control problem and solves it using tools from
optimal control theory \cite{doi:10.1137/19M1288802}. Recent
neural-network-based methods further provide flexible tools for solving such
optimal control formulations. In particular, Wei et al. developed a
neural-network-based numerical method for the nondegenerate case
\cite{wei_optimal_2022}, while Chen et al. proposed a method combining
successive approximations with nested neural networks based on the Pontryagin
maximum principle \cite{chen_detecting_2024}. Inspired by these developments,
we extend the optimal-control-based neural-network approach to degenerate
stochastic systems with multiplicative Gaussian noise and develop a numerical
method for the resulting degenerate optimal control problem.

The main contributions of this paper are threefold. First, we derive an OM
functional for degenerate stochastic dynamical systems with multiplicative
Gaussian noise, a setting in which the diffusion coefficient is both degenerate
and state dependent. The derivation first treats the time-dependent
additive-noise case via the Girsanov transformation and then uses a suitable
Lamperti-type transformation \cite{moller2010state} to handle the state dependence of the diffusion
coefficient. Second, using the derived OM functional, we formulate the MPTP
problem as a constrained variational problem and derive the corresponding
Hamilton--Pontryagin boundary value problem. This extends the degenerate
variational characterization from the additive Gaussian noise case
\cite{WOS:001514271800001} to the multiplicative Gaussian noise case. Third, we
reformulate this constrained variational problem as an equivalent optimal
control problem and develop a neural-network-based numerical method for
computing MPTPs in degenerate stochastic dynamical systems with multiplicative
Gaussian noise.

The remainder of this paper is organized as follows. Section~\ref{sec2} collects the preliminaries, including the model formulation, the definition of the OM functional, the standing assumptions, and several technical lemmas. Section~\ref{sec3} presents the main theoretical results: the derivation of the OM functional, the Hamilton--Pontryagin characterization of the MPTP problem, the optimal control reformulation, and the proofs of the main theorems. Section~\ref{sec4} provides several examples illustrating the theoretical results. Section~\ref{sec5} concludes the paper. Additional technical details are deferred to the Appendix.

\section{Preliminaries}
\label{sec2}

In this section, we introduce the degenerate stochastic system studied in this paper and recall the notions of the OM functional and the MPTP. We then state the standing assumptions and collect the technical lemmas needed in the subsequent analysis.

\subsection{Degenerate stochastic systems with multiplicative Gaussian noise}

We first fix the notation for the regularity class used below. For
\(r\in\mathbb N\), we write
\[
h\in C_b^{0,r}
\bigl([0,T]\times\mathbb R^d\times\mathbb R^m;\mathbb R^k\bigr)
\]
if \(h\) is continuous in \(t\), \(r\)-times continuously differentiable in the
spatial variables \((x,y)\), and \(h\), together with all its spatial
derivatives up to order \(r\), is bounded on
\([0,T]\times\mathbb R^d\times\mathbb R^m\).

%

Let \((\Omega, \mathcal{F}, \{\mathcal{F}_t\}_{t \geq 0}, \mathbb{P})\) be a complete filtered probability space, where \(\{\mathcal{F}_t\}_{t \geq 0}\) is a nondecreasing family of sub-\(\sigma\)-fields of \(\mathcal{F}\) satisfying the usual conditions. We consider the following degenerate stochastic system with multiplicative Gaussian noise on \([0,T]\):
\begin{equation}\label{E-1}
    \begin{cases}
   \mathrm{d}X_{t}=g(t,X_{t},Y_{t})\,\mathrm{d}t,\qquad X_{0}=x_0 \in \mathbb{R}^{d},\\
   \mathrm{d}Y_{t}=f(t,X_{t},Y_{t})\,\mathrm{d}t+\sigma(t,X_t,Y_t)\,\mathrm{d}W_{t},\qquad Y_{0}=y_0 \in \mathbb{R}^{m}.
    \end{cases}
\end{equation}
Here \(W_t\) is an \(m\)-dimensional standard Brownian motion. The coefficients
are assumed to satisfy
\(
g\in C_b^{0,1}
\bigl([0,T]\times\mathbb R^d\times\mathbb R^m;\mathbb R^d\bigr),
f\in C_b^{0,2}
\bigl([0,T]\times\mathbb R^d\times\mathbb R^m;\mathbb R^m\bigr),
\)
and
\(
\sigma,\sigma^{-1}\in C_b^{0,2}
\bigl([0,T]\times\mathbb R^d\times\mathbb R^m;\mathbb R^{m\times m}\bigr).
\)
Here \(\sigma^{-1}\) denotes the matrix inverse of \(\sigma\). In particular,
the diffusion matrix is uniformly nondegenerate in the \(Y\)-component: there
exists a constant \(\lambda>0\) such that
\[
\xi^\top \sigma(t,x,y)\sigma(t,x,y)^\top \xi
\geq \lambda |\xi|^2,
\qquad
\text{for all } (t,x,y)\in[0,T]\times\mathbb R^d\times\mathbb R^m
\text{ and } \xi\in\mathbb R^m.
\]
Nevertheless, system~\eqref{E-1} is degenerate as an \((X,Y)\)-system because
the Brownian noise acts only on the \(Y\)-equation and no Brownian component
appears directly in the \(X\)-equation. Under these assumptions,
system~\eqref{E-1} admits a unique strong solution; see Lan and
Wu~\cite{lan_new_2014}.

%

\subsection{Onsager--Machlup functional and most probable transition pathways}
\label{sec2.2}

For notational convenience, set
\(
Z_t=(X_t,Y_t)^\top, z_0=(x_0,y_0)^\top\in\mathbb{R}^{d+m}.
\)
Let \(C^{z_0}([0,T];\mathbb{R}^{d+m})\) denote the space of continuous paths starting from \(z_0\), equipped with the uniform norm
\[
\|z\|:=\|z\|_{[0,T]}=\sup_{t\in[0,T]} |z(t)|.
\]
We also define the Cameron--Martin space
\[
\mathcal{H}
=
\left\{
h:[0,T]\to\mathbb{R}^{d+m}:
h \text{ is absolutely continuous},\
\dot h\in L^2([0,T];\mathbb{R}^{d+m}),\
h(0)=0
\right\}.
\]
For \(\phi\in \mathcal{H}_{z_0}:=z_0+\mathcal{H}\) and \(\varepsilon>0\), define
the open tube around \(\phi\) by
\[
K(\phi,\varepsilon)
=
\left\{
z\in C^{z_0}([0,T];\mathbb{R}^{d+m}):
\|z-\phi\|< \varepsilon
\right\}.
\]
\begin{definition}[Onsager--Machlup function and action functional]
Let \(Z\) be the solution process of system \eqref{E-1}. If, for sufficiently small
\(\varepsilon>0\), the tube probability around a reference path \(\phi\in \mathcal{H}_{z_0}\)
admits the asymptotic form
\begin{equation*}\label{defom}
\mathbb{P}\bigl(Z\in K(\phi,\varepsilon)\bigr)
\propto
C(\varepsilon)
\exp\left\{
-\int_0^T
\operatorname{OM}(\phi,\dot{\phi})\,\mathrm{d}t
\right\},
\qquad \varepsilon\to0^+,
\end{equation*}
where \(C(\varepsilon)\) is independent of \(\phi\), then
\(\operatorname{OM}(\phi,\dot{\phi})\) is called the Onsager--Machlup function. The functional
\(
I(\phi)
=
\int_0^T
\operatorname{OM}(\phi,\dot{\phi})\,\mathrm{d}t
\)
is called the corresponding Onsager--Machlup action functional.
\end{definition}

For the degenerate systems considered here, admissible transition paths must satisfy the constraint imposed by the unforced component. Writing
\[
\phi(t)=(\phi_1(t),\phi_2(t))^\top,
\qquad
\phi_1(t)\in\mathbb{R}^d,\quad
\phi_2(t)\in\mathbb{R}^m,
\]
we define the admissible class connecting \(z_0\) to a prescribed terminal point \(z_T\) by
\[
\mathcal{A}_{z_0,z_T}
=
\left\{
\phi\in\mathcal{H}_{z_0}:
\phi(T)=z_T,\ 
\dot{\phi}_1(t)=g(t,\phi_1(t),\phi_2(t))
\text{ for a.e. } t\in[0,T]
\right\}.
\]

\begin{definition}[Most probable transition pathway]
\normalfont
A path \(\phi^*\in\mathcal{A}_{z_0,z_T}\) is called a most probable transition pathway from \(z_0\) to \(z_T\) if
\[
I(\phi^*)
=
\inf_{\phi\in\mathcal{A}_{z_0,z_T}} I(\phi).
\]
\end{definition}

Thus, in the present degenerate setting, the MPTP problem is naturally formulated as a constrained variational problem. The explicit Onsager--Machlup functional and the associated variational equations will be derived in the next section.

\subsection{Assumptions}

To characterize the Onsager--Machlup functional through small-tube asymptotics, we impose the following assumptions, adapted from Liu and Gao~\cite{liu_onsagermachlup_2024}. Throughout this subsection, \(W=(W_1,\ldots,W_m)\) denotes an \(m\)-dimensional Brownian motion.

\begin{assumption}
\normalfont
\label{ass:H1}
For every \(1\leq i\leq m\),
\[
\|(W_1,\ldots,W_i,\ldots,W_m)\|
=
\|(W_1,\ldots,-W_i,\ldots,W_m)\|.
\]
\end{assumption}

\begin{assumption}
\normalfont
\label{ass:H2}
For every \(1\leq i\leq m\) and every \(c\in\mathbb{R}\),
\[
\limsup_{\varepsilon\to0}
\mathbb{E}\left[
\exp\bigl(c|W_i(T)|^2\bigr)
\,\middle|\, \|W\|< \varepsilon
\right]
\leq 1.
\]
\end{assumption}

\begin{assumption}
\normalfont
\label{ass:H3}
There exist constants \(p>0\), \(q\in(0,p\wedge4)\), and \(C_1,C_2>0\) such that, for every sufficiently small \(\varepsilon>0\), the following two conditions hold:
\[
\|W\|< \varepsilon
\quad\Longrightarrow\quad
\int_0^T |W(t)|^4\,\mathrm{d}t
\leq C_1\varepsilon^{p\wedge 4},
\]
and
\[
\mathbb{P}(\|W\|< \varepsilon)
\geq
\exp\left(-\frac{C_2}{\varepsilon^q}\right).
\]
\end{assumption}


\subsection{Technical lemmas}

We finally recall three technical lemmas that will be used in the derivation of the Onsager--Machlup functional.

\begin{lemma}[Ikeda and Watanabe \cite{Ikeda2014stochastic}]
\normalfont
\label{lem:ikeda_watanabe}
For a fixed \(n\geq1\), let \(I_1,\ldots,I_n\) be random variables defined on \((\Omega,\mathcal{F},\mathbb{P})\), and let \(\{A_\varepsilon:\varepsilon>0\}\) be a family of events in \(\mathcal{F}\). Suppose that, for every \(c\in\mathbb{R}\) and every \(i=1,\ldots,n\),
\[
\limsup_{\varepsilon\to0}
\mathbb{E}\left[
\exp(cI_i)\,\middle|\, A_\varepsilon
\right]
\leq 1.
\]
Then
\[
\lim_{\varepsilon\to0}
\mathbb{E}\left[
\exp\left(\sum_{i=1}^n I_i\right)
\,\middle|\, A_\varepsilon
\right]
=1.
\]
\end{lemma}

\begin{lemma}[Borell~\cite{borell1977gaussian}; Ledoux~\cite{ledoux1996isoperimetry}]
\normalfont
\label{lem:borell_ledoux}
Let \(f\in L^2([0,1])\) be deterministic and define
\(
I_i(f)=\int_0^1f(t)\,\mathrm{d}W_i(t).
\)
Then
\[
\lim_{\varepsilon\to0}
\mathbb{E}\left[
\exp\bigl(|I_i(f)|\bigr)
\,\middle|\, \|W\|<\varepsilon
\right]
=1.
\]
\end{lemma}

\begin{lemma}[Shepp and Zeitouni~\cite{sheppzeitouni1993}]
\normalfont
\label{lem:shepp_zeitouni}
Let \(1\leq i\leq m\). Let \(\mathcal{F}_i\) be the \(\sigma\)-field generated by
\[
\bigl\{
W_1(t),\ldots,W_{i-1}(t),W_{i+1}(t),\ldots,W_m(t):
0\leq t\leq 1
\bigr\}.
\]
Suppose that \(\Psi\) is an \(\mathcal{F}_i\)-adapted process satisfying, for
every \(c>0\),
\[
\lim_{\varepsilon\to0}
\mathbb{E}\left[
\exp\left\{
c\int_0^1 \Psi^2(t)\,\mathrm{d}t
\right\}
\,\middle|\, \|W\|<\varepsilon
\right]
=1.
\]
Then
\[
\lim_{\varepsilon\to0}
\mathbb{E}\left[
\exp\left\{
\left|\int_0^1 \Psi(t)\,\mathrm{d}W_i(t)\right|
\right\}
\,\middle|\, \|W\|<\varepsilon
\right]
=1.
\]
\end{lemma}

\begin{remark}
\label{rem:time_rescaling}
By the Brownian scaling property, if
\(
\widetilde W(t):=T^{-1/2}W(Tt), 0\leq t\leq1,
\)
then \(\widetilde W\) is again a Brownian motion on \([0,1]\), and
\[
\|W\|_{[0,T]}
=
\sqrt T\,\|\widetilde W\|_{[0,1]}.
\]
The same time-rescaling transforms stochastic integrals on \([0,T]\) into
stochastic integrals on \([0,1]\). Hence Lemmas~\ref{lem:borell_ledoux} and
\ref{lem:shepp_zeitouni} also hold on any finite interval \([0,T]\), with the
supremum norm taken over \([0,T]\).
\end{remark}


\section{General Theory}\label{sec3}

\subsection{Onsager--Machlup theory for degenerate stochastic systems with multiplicative Gaussian noise}\label{sec3.1}

To derive the Onsager--Machlup functional for system~\eqref{E-1}, we first consider the following auxiliary system with additive noise:
\begin{equation}\label{E-2}
    \begin{cases}
   \mathrm{d}X_{t}=g(t,X_{t},Y_{t})\,\mathrm{d}t,\qquad
   X_{0}=x_0 \in \mathbb{R}^{d},\\
   \mathrm{d}Y_{t}=f(t,X_{t},Y_{t})\,\mathrm{d}t
   +\sigma(t)\,\mathrm{d}W_{t},\qquad
   Y_{0}=y_0 \in \mathbb{R}^{m},
    \end{cases}
\end{equation}
where \(\sigma\in C_b^2([0,T];\mathbb{R}^{m\times m})\) is uniformly nondegenerate. The Onsager--Machlup functional is first derived for system~\eqref{E-2}, and then extended to system~\eqref{E-1} by means of a
suitable transformation.

We now state the main results. The proofs are given in Section~\ref{sec3.4}.


\begin{theorem}\label{theorem 3.1}
Consider the stochastic dynamical system~\eqref{E-2}. Assume that
 \(
g\in C_b^{0,1}
\bigl([0,T]\times\mathbb R^d\times\mathbb R^m;\mathbb R^d\bigr),
f\in C_b^{0,2}
\bigl([0,T]\times\mathbb R^d\times\mathbb R^m;\mathbb R^m\bigr).
\)
Assume further that \(\sigma:[0,T]\to\mathbb R^{m\times m}\) is \(C^2\) with bounded derivatives up to order two and is uniformly nondegenerate. Then its Onsager--Machlup action functional is
    \begin{equation}
        I(\phi)
        =
        \int_0^T \operatorname{OM}(\phi,\dot{\phi})\,\mathrm{d}t,
    \end{equation}
    where
    \begin{equation}\label{OM1}
        \operatorname{OM}(\phi,\dot{\phi})
        =
        \frac{1}{2}
        \Big\lvert
        \sigma^{-1}(t)
        \big(\dot{\phi}_2(t)-f(t,\phi(t))\big)
        \Big\rvert^2
        +
        \frac{1}{2}\operatorname{div}_y f(t,\phi(t)).
    \end{equation}
    Here \(\phi=(\phi_1,\phi_2)^T\) is an absolutely continuous reference path satisfying \(\phi(0)=z_0=(x_0,y_0)^T\) and the constraint
    \(
        \dot{\phi}_1(t)=g(t,\phi_1(t),\phi_2(t)).
    \)
\end{theorem}

\begin{remark}\label{re3.1}
When \(\sigma(t)\) is constant, Formula~\eqref{OM1} reduces to the usual Onsager--Machlup functional for additive-noise degenerate systems with a constant diffusion matrix \cite{capitaine_onsager-machlup_1995,liu_onsagermachlup_2024}. Thus, the present formulation extends the known additive-noise result to the case where the diffusion matrix depends explicitly on time, namely \(\sigma=\sigma(t)\).
\end{remark}

Before stating the result for the multiplicative-noise system~\eqref{E-1}, we introduce a Lamperti-type transform \cite{moller2010state}. Since such a transform may only be
defined locally, we work on an open domain containing the reference path.
\begin{assumption}[Lamperti domain]\label{ass:lamperti-domain}
There exists an open set
\(
    \mathcal D\subset \mathbb R^d\times\mathbb R^m
\)
such that, on \([0,T]\times\mathcal D\), the diffusion matrix \(\sigma\) is
nondegenerate, that is,
\[
    \det \sigma(t,x,y)\neq 0,
    \qquad
    (t,x,y)\in[0,T]\times\mathcal D,
\]
and
\(
    \sigma^{-1}\in C_b^{0,2}
    \bigl([0,T]\times\mathcal D;\mathbb R^{m\times m}\bigr).
\)
Moreover, there exists a map
\(
    \Gamma:[0,T]\times\mathcal D\to\mathbb R^m
\)
which is sufficiently smooth, with the derivatives appearing below bounded
on \([0,T]\times\mathcal D\), such that
\begin{equation}\label{condition1}
    \partial_y \Gamma(t,x,y)=\sigma^{-1}(t,x,y),
    \qquad
    (t,x,y)\in[0,T]\times\mathcal D.
\end{equation}
For each \(t\in[0,T]\), the map
\(
    \Theta_t:\mathcal D\to \Theta_t(\mathcal D),
    \Theta_t(x,y)=(x,\Gamma(t,x,y)),
\)
is assumed to be a \(C^2\)-diffeomorphism onto its image.
\end{assumption}

\begin{theorem}\label{theorem3.2}
Consider system~\eqref{E-1}. Assume that
\(
g\in C_b^{0,1}
\bigl([0,T]\times\mathbb R^d\times\mathbb R^m;\mathbb R^d\bigr),
f\in C_b^{0,2}
\bigl([0,T]\times\mathbb R^d\times\mathbb R^m;\mathbb R^m\bigr),
\)
and
\(
\sigma\in C_b^{0,2}
\bigl([0,T]\times\mathbb R^d\times\mathbb R^m;\mathbb R^{m\times m}\bigr).
\)
Suppose that Assumption~\ref{ass:lamperti-domain} holds on an open Lamperti
domain \(\mathcal D\). Let
\(
    \phi(t)=(\phi_1(t),\phi_2(t))^T,
t\in[0,T],
\)
be an admissible reference path satisfying
\(
    \phi(0)=z_0=(x_0,y_0)^T,
\dot{\phi}_1(t)=g(t,\phi_1(t),\phi_2(t)).
\)
Assume further that
\(
    \phi([0,T])\Subset\mathcal D,
\)
where \(\Subset\) denotes compact containment in \(\mathcal D\). Then the
Onsager--Machlup function of system~\eqref{E-1} along the path \(\phi\) is
given by
\begin{equation}\label{OM2}
    \operatorname{OM}(\phi,\dot{\phi})
    =
    \frac12
    \left|
        \sigma^{-1}
        \left[
            \dot{\phi}_2
            -
            f
            -
            \frac12
            \sigma
            \operatorname{Tr}
            \left(
                \partial_{yy}^2\Gamma\,\sigma\sigma^T
            \right)
        \right]
    \right|^2
    +
    \frac12
    \operatorname{tr}
    \left(
        \partial_y H\sigma
    \right),
\end{equation}
where all the coefficients in~\eqref{OM2} are evaluated at
\((t,\phi_1(t),\phi_2(t))\). The transformed drift \(H\) is defined by
\begin{align}\label{H-def}
    H(t,x,y)
    =
    &\,
    \partial_t \Gamma(t,x,y)
    +
    \partial_x \Gamma(t,x,y)g(t,x,y)
    +
    \sigma^{-1}(t,x,y)f(t,x,y) \notag \\
    &+
    \frac12
    \operatorname{Tr}
    \left(
        \partial_{yy}^2\Gamma(t,x,y)
        \sigma\sigma^T(t,x,y)
    \right).
\end{align}

In the above formulas,
\(\partial_x \Gamma=(\partial_{x_j}\Gamma_i)_{1\le i\le m,\,1\le j\le d}\) and
\(\partial_y \Gamma=(\partial_{y_j}\Gamma_i)_{1\le i,j\le m}\) denote the Jacobian
matrices with respect to \(x\) and \(y\), respectively. Moreover,
\(\partial_{yy}^2\Gamma_i\) denotes the Hessian matrix of the \(i\)-th component
\(\Gamma_i\) with respect to \(y\), and
\[
    \operatorname{Tr}
    \left(
        \partial_{yy}^2\Gamma\,\sigma\sigma^T
    \right)
    :=
    \left(
        \operatorname{tr}
        \left(
            \partial_{yy}^2\Gamma_1\,\sigma\sigma^T
        \right),
        \ldots,
        \operatorname{tr}
        \left(
            \partial_{yy}^2\Gamma_m\,\sigma\sigma^T
        \right)
    \right)^T .
\]
\end{theorem}

\begin{remark}\label{re-lamperti}
Assumption~\ref{ass:lamperti-domain} means that, at least on a neighborhood of
the reference path, the multiplicative noise can be converted into additive
noise by the Lamperti-type change of variables
\[
    (x,y)\mapsto (x,\Gamma(t,x,y)).
\]
Indeed, by \(\partial_y\Gamma=\sigma^{-1}\), the diffusion coefficient in the
transformed variable becomes the identity matrix. The drift \(H\) in
\eqref{H-def} is precisely the drift term obtained from It\^o's formula after
this change of variables. Hence \eqref{OM2} is obtained by applying the
additive-noise Onsager--Machlup formula in Theorem~\ref{theorem 3.1} to the
transformed system and then rewriting the result in the original variables.
\end{remark}

\begin{remark}\label{re-scalar-case}
The formula~\eqref{OM2} is formulated under the existence of a local
Lamperti-type transform \(\Gamma\) satisfying \eqref{condition1}. In special cases
where \(\Gamma\) can be written explicitly, \eqref{OM2} reduces to the familiar
Onsager--Machlup functionals.

For example, when \(d=m=1\), the identity \(\Gamma_y=\sigma^{-1}\) gives
\(
    \Gamma_{yy}
    =
    -\frac{\sigma_y}{\sigma^2}.
\)
Thus, with all coefficients evaluated at \((t,\phi_1(t),\phi_2(t))\),
\eqref{OM2} becomes
\begin{equation}\label{OM3}
    \operatorname{OM}(\phi,\dot{\phi})
    =
    \frac{1}{2}
    \left(
    \frac{\dot{\phi}_2-f+\frac{1}{2}\sigma\sigma_y}{\sigma}
    \right)^2
    +
    \frac{1}{2}\sigma
    \partial_y
    \left(
    \Gamma_t+g \Gamma_x+\frac{f}{\sigma}
    -\frac{1}{2}\sigma_y
    \right).
\end{equation}
Here \(\Gamma_t=\partial_t\Gamma\), \(\Gamma_x=\partial_x\Gamma\), and
\(\sigma_y=\partial_y\sigma\).

If, in addition, \(\sigma=\sigma(y)\), then one may choose, locally,
\(
    \Gamma(y)=\int^y \frac{1}{\sigma(r)}\,\mathrm dr .
\)
In this case \(\Gamma_t=\Gamma_x=0\), and \eqref{OM3} coincides with the classical
one-dimensional multiplicative-noise Onsager--Machlup formula, such as
equation~(34) in \cite{horsthemke_onsager-machlup_1975}.
\end{remark}

\subsection{Characterization of most probable transition pathways}\label{sec3.2}

Throughout this subsection, the transition time \(T>0\) is fixed.
Given two endpoints \(z_0=(x_0,y_0)^T\) and \(z_T=(x_T,y_T)^T\) in the Lamperti
domain \(\mathcal D\), we seek an MPTP
\(\phi=(\phi_1,\phi_2)^T\) connecting \(z_0\) to \(z_T\) over the prescribed
time interval \([0,T]\), that is,
\[
    \phi(0)=z_0,
    \quad
    \phi(T)=z_T,
    \quad
    \phi([0,T])\Subset\mathcal D.
\]
As discussed in Section~\ref{sec2.2}, an MPTP can be characterized as a minimizer of the Onsager--Machlup action functional. In the nondegenerate case, the Onsager--Machlup function can be regarded as a classical Lagrangian, and the corresponding Euler--Lagrange equation gives a two-point boundary value problem. For the degenerate system considered in this paper, however, the \(X\)-component contains no direct noise. Therefore,
the constraint
\(
    \dot{\phi}_1(t)=g(t,\phi_1(t),\phi_2(t))
\)
must be imposed explicitly. Following the approach developed in our previous work on degenerate additive-noise systems~\cite{WOS:001514271800001}, this constrained variational problem can be treated by introducing a Lagrange multiplier and applying the
Hamilton--Pontryagin variational principle.


For system~\eqref{E-1}, we take the Onsager--Machlup function in \eqref{OM2} as the Lagrangian and denote it by
\[
    \mathcal{L}_{\mathrm{OM}}(t,\phi_1,\phi_2,\dot{\phi}_2)
    :=
    \operatorname{OM}(\phi,\dot{\phi}).
\]
Explicitly, it is given by
\begin{align*}\label{LOM}
    \mathcal{L}_{\mathrm{OM}}
    (t,\phi_1,\phi_2,\dot{\phi}_2)
    &=
    \frac{1}{2}
    \left|
    \sigma^{-1}
    \left[
    \dot{\phi}_2
    -
    f
    -
    \frac{1}{2}\sigma
    \operatorname{Tr}
    \left(
    \partial_{yy}^2\Gamma\,\sigma\sigma^T
    \right)
    \right]
    \right|^2
    +
    \frac12
    \operatorname{tr}
    \left(
        \partial_y H\sigma
    \right),
\end{align*}
where all coefficients are evaluated at \((t,\phi(t))\).
The constrained variational problem associated with the MPTP is to find a
stationary curve of the action
\begin{equation*}\label{CV}
    I(\phi)
    =
    \int_0^T
    \mathcal{L}_{\mathrm{OM}}
    (t,\phi_1,\phi_2,\dot{\phi}_2)
    \,\mathrm{d}t
\end{equation*}
over all admissible paths satisfying
\begin{equation*}\label{CV-constraints}
\dot{\phi}_1(t)=g(t,\phi_1(t),\phi_2(t)),
    \quad
    \phi(0)=z_0,
    \quad
    \phi(T)=z_T,
    \quad
    \phi([0,T])\Subset\mathcal D .
\end{equation*}

To treat this constrained variational problem, we introduce a time-dependent Lagrange multiplier \(\lambda:[0,T]\to\mathbb{R}^d\) associated with the constraint \(\dot{\phi}_1=g(t,\phi_1,\phi_2)\). This leads to the
Hamilton--Pontryagin variational principle.

\begin{pro}[Hamilton--Pontryagin principle]\label{prop:HP}
Let
\(
    \gamma=(\phi_1,\dot{\phi}_1,\phi_2,\dot{\phi}_2,\lambda)
\)
be a smooth curve with \(\lambda:[0,T]\to\mathbb{R}^d\) , satisfying
\(
    \phi(0)=z_0,
    \phi(T)=z_T,
  \phi([0,T])\Subset\mathcal D.
\)
The critical point condition \(\delta\mathcal{J}=0\) for the
Hamilton--Pontryagin functional
\begin{equation}\label{HP-functional}
\mathcal{J}(\phi_1,\dot{\phi}_1,\phi_2,\dot{\phi}_2,\lambda)
=
\int_0^T
\mathcal{L}_{\mathrm{OM}}
(t,\phi_1,\phi_2,\dot{\phi}_2)\,\mathrm{d}t 
+
\int_0^T
\left\langle
\lambda,
g(t,\phi_1,\phi_2)-\dot{\phi}_1
\right\rangle
\,\mathrm{d}t
\end{equation}
yields the Hamilton--Pontryagin system
\begin{equation}\label{HP}
\left\{
\begin{aligned}
-
\frac{\mathrm{d}}{\mathrm{d}t}
\bigl(
\nabla_{\dot{\phi}_2}\mathcal{L}_{\mathrm{OM}}
\bigr)
+
\nabla_{\phi_2}\mathcal{L}_{\mathrm{OM}}
+
(D_{\phi_2}g)^T\lambda
&=0,\\
\nabla_{\phi_1}\mathcal{L}_{\mathrm{OM}}
+
(D_{\phi_1}g)^T\lambda
+
\dot{\lambda}
&=0,\\
g(t,\phi_1,\phi_2)-\dot{\phi}_1
&=0.
\end{aligned}
\right.
\end{equation}
Here all derivatives of \(\mathcal{L}_{\mathrm{OM}}\) are evaluated at \((t,\phi_1(t),\phi_2(t),\dot{\phi}_2(t))\), whereas \(D_{\phi_1}g\) and
\(D_{\phi_2}g\) are evaluated at \((t,\phi_1(t),\phi_2(t))\). The gradients
\(
\nabla_{\phi_1}\mathcal{L}_{\mathrm{OM}}\in\mathbb{R}^d,
\nabla_{\phi_2}\mathcal{L}_{\mathrm{OM}}\in\mathbb{R}^m,
\nabla_{\dot{\phi}_2}\mathcal{L}_{\mathrm{OM}}\in\mathbb{R}^m
\)
are regarded as column vectors, and
\(
D_{\phi_1}g\in\mathbb{R}^{d\times d},
D_{\phi_2}g\in\mathbb{R}^{d\times m}
\)
denote the Jacobian matrices of \(g\) with respect to \(\phi_1\) and \(\phi_2\), respectively.
\end{pro}

\begin{proof}
Taking variations of \(\mathcal{J}\) with respect to \(\lambda\) gives
\[
    g(t,\phi_1,\phi_2)-\dot{\phi}_1=0.
\]
Next, taking variations with respect to \(\phi_1\) and integrating by parts,
with \(\delta\phi_1(0)=\delta\phi_1(T)=0\), yields
\[
    \nabla_{\phi_1}\mathcal{L}_{\mathrm{OM}}
    +
    (D_{\phi_1}g)^T\lambda
    +
    \dot{\lambda}
    =0.
\]
Similarly, variations with respect to \(\phi_2\), together with
\(\delta\phi_2(0)=\delta\phi_2(T)=0\), yield
\[
    -
    \frac{\mathrm{d}}{\mathrm{d}t}
    \nabla_{\dot{\phi}_2}\mathcal{L}_{\mathrm{OM}}
    +
    \nabla_{\phi_2}\mathcal{L}_{\mathrm{OM}}
    +
    (D_{\phi_2}g)^T\lambda
    =0.
\]
Thus \((\phi,\lambda)\) satisfies \eqref{HP}. This proves the result.
\end{proof}

\subsection{Optimal control formulation of the MPTP}\label{sec3.3}

Although Proposition~\ref{prop:HP} characterizes stationary curves through
the Hamilton--Pontryagin boundary value problem~\eqref{HP}, directly solving
the resulting coupled system can be computationally challenging, especially
in high-dimensional settings. To facilitate numerical computation, we
reformulate the minimization of the Onsager--Machlup action as a constrained
optimal control problem. In this formulation, the control represents the
noise realization that drives the system from the initial state to the
prescribed terminal state over the fixed time interval \([0,T]\).

For comparison, consider a nondegenerate stochastic system in \(\mathbb{R}^r\)
\[
    \mathrm{d}P_t
    =
    b(t,P_t)\,\mathrm{d}t
    +
    \Sigma(t,P_t)\,\mathrm{d}W_t .
\]
The corresponding controlled dynamics are
\[
    \dot{\phi}(t)
    =
    b(t,\phi(t))
    +
    \Sigma(t,\phi(t))\theta(t),
    \qquad
    \theta\in L^2([0,T];\mathbb{R}^r).
\]
In this case, the Onsager--Machlup action may be written in terms of the
control by defining
\[
    \ell_{\mathrm{OM}}^{\mathrm{nd}}(t,\phi,\theta)
    =
    \mathcal{L}_{\mathrm{OM}}
    \bigl(
    t,\phi,
    b(t,\phi)+\Sigma(t,\phi)\theta
    \bigr).
\]
The corresponding constrained optimal control problem is then
\begin{equation}\label{OC}
\left\{
\begin{aligned}
\inf_{\phi,\theta}\quad
&
J(\phi,\theta)
=
\int_0^T
\ell_{\mathrm{OM}}^{\mathrm{nd}}(t,\phi(t),\theta(t))
\,\mathrm{d}t
+
\eta(\phi(T)),
\\
\text{subject to}\quad
&
\dot{\phi}(t)
=
b(t,\phi(t))
+
\Sigma(t,\phi(t))\theta(t),
\\
&
\phi(0)=x_0.
\end{aligned}
\right.
\end{equation}
Here \(\eta\) enforces the prescribed terminal condition. For example, for a
given terminal point \(x_T\), one may take
\[
    \eta(x)
    =
    \begin{cases}
    0, & x=x_T,\\
    +\infty, & x\neq x_T.
    \end{cases}
\]
In numerical computations, this hard terminal constraint is often replaced
by a penalty term such as \(\tau_T|\phi(T)-x_T|^2\), where \(\tau_T>0\) is
chosen sufficiently large.

For the degenerate multiplicative noise system~\eqref{E-1}, we introduce a
control \(\theta\in L^2([0,T];\mathbb{R}^m)\) acting only in the forced
component. The controlled dynamics are
\[
\left\{
\begin{aligned}
    \dot{\phi}_1(t)
    &=
    g(t,\phi_1(t),\phi_2(t)),
    \\
    \dot{\phi}_2(t)
    &=
    f(t,\phi_1(t),\phi_2(t))
    +
    \sigma(t,\phi_1(t),\phi_2(t))\theta(t).
\end{aligned}
\right.
\]
The admissible state trajectories are required to remain in the Lamperti
domain; more precisely, we impose
\(
    \phi([0,T])\Subset\mathcal D .
\)
We assume that the endpoint \(z_T\) is reachable under these controlled
dynamics by an admissible trajectory satisfying
\(\phi([0,T])\Subset\mathcal D\); otherwise, the terminal condition is
understood in the penalized sense. The corresponding constrained optimal
control formulation of the OM action minimization problem on the fixed time
interval \([0,T]\) is
\begin{equation}\label{OC-degenerate}
\left\{
\begin{aligned}
\inf_{\phi,\theta}\quad
&
J(\phi,\theta)
=
\int_0^T
\ell_{\mathrm{OM}}(t,\phi(t),\theta(t))
\,\mathrm{d}t
+
\eta(\phi(T)),
\\
\text{subject to}\quad
&
\dot{\phi}_1(t)
=
g(t,\phi_1(t),\phi_2(t)),
\\
&
\dot{\phi}_2(t)
=
f(t,\phi_1(t),\phi_2(t))
+
\sigma(t,\phi_1(t),\phi_2(t))\theta(t),
\\
&
\phi(0)=z_0,
\\
&
\phi([0,T])\Subset\mathcal D.
\end{aligned}
\right.
\end{equation}
Here \(\eta\) is chosen as the indicator function of \(\{z_T\}\) for the hard
terminal constraint \(\phi(T)=z_T\), or is replaced by a large terminal
penalty in numerical implementations.

The running cost \(\ell_{\mathrm{OM}}\) is obtained by substituting the
controlled velocity
\(
    \dot{\phi}_2
    =
    f(t,\phi)
    +
    \sigma(t,\phi)\theta
\)
into the Onsager--Machlup Lagrangian. Namely,
\[
    \ell_{\mathrm{OM}}(t,\phi,\theta)
    =
    \mathcal{L}_{\mathrm{OM}}
    \bigl(
    t,\phi_1,\phi_2,
    f(t,\phi)+\sigma(t,\phi)\theta
    \bigr).
\]
Therefore,
\begin{align}\label{ellOM}
    \ell_{\mathrm{OM}}(t,\phi,\theta)
    =
    \frac{1}{2}
    \left|
    \theta
    -
    \frac{1}{2}
    \operatorname{Tr}
    \left(
    \partial_{yy}^2\Gamma\,\sigma\sigma^T
    \right)
    \right|^2
    +
    \frac12
    \operatorname{tr}
    \left(
        \partial_y H\sigma
    \right),
\end{align}
where all coefficients are evaluated at
\((t,\phi(t))\in[0,T]\times\mathcal D\), and the notation \(u\), \(H\), and
\(\operatorname{div}_u H\) is the same as in the definition of
\(\mathcal{L}_{\mathrm{OM}}\). The trace term is understood componentwise,
that is,
\[
    \left[
    \operatorname{Tr}
    \left(
    \partial_{yy}^2\Gamma\,\sigma\sigma^T
    \right)
    \right]_k
    =
    \operatorname{Tr}
    \left(
    \partial_{yy}^2\Gamma_k\,\sigma\sigma^T
    \right),
    \qquad k=1,\ldots,m.
\]
This formulation should be distinguished from the classical
Freidlin--Wentzell optimal control problem~\cite{wei_optimal_2022}, whose
running cost is \(\frac{1}{2}|\theta|^2\). Here the running cost is the
Onsager--Machlup density after substituting the controlled dynamics.

If an optimal state trajectory \(\phi^*\) of~\eqref{OC-degenerate} exists, then it is identified with a minimizer of the OM action over admissible paths
satisfying \(\phi([0,T])\Subset\mathcal D\), and hence with an MPTP connecting \(z_0\) to \(z_T\) over the prescribed time interval \([0,T]\).

\paragraph{Neural-network approximation.}
In practice, both the path \(\phi\) and the control \(\theta\) may be approximated
by neural networks, denoted by \(\phi^{NN}\) and \(\theta^{NN}\), respectively.
Let \(\{t_i\}_{i=1}^{N_T}\subset[0,T]\) be collocation points. In implementations based on the local Lamperti
representation, the loss is evaluated for approximated trajectories whose
sampled values remain in \(\mathcal D\); if necessary, an additional barrier
penalty can be introduced to discourage exits from \(\mathcal D\). A typical
loss function for~\eqref{OC-degenerate} is
\begin{align}\label{NN-loss}
    \operatorname{Loss}
    =
    \tau_{\phi}\operatorname{Loss}_{\phi}
    +
    \tau_{\mathrm{OM}}\operatorname{Loss}_{\mathrm{OM}}
    +
    \tau_b\operatorname{Loss}_b,
\end{align}
where
\begin{align*}
    \operatorname{Loss}_{\phi}
    &=
    \operatorname{Loss}_{\phi_1}
    +
    \operatorname{Loss}_{\phi_2},
    \\
    \operatorname{Loss}_{\phi_1}
    &=
    \frac{1}{N_T}
    \sum_{i=1}^{N_T}
    \left|
    \dot{\phi}^{NN}_1(t_i)
    -
    g(t_i,\phi^{NN}(t_i))
    \right|^2,
    \\
    \operatorname{Loss}_{\phi_2}
    &=
    \frac{1}{N_T}
    \sum_{i=1}^{N_T}
    \left|
    \dot{\phi}^{NN}_2(t_i)
    -
    f(t_i,\phi^{NN}(t_i))
    -
    \sigma(t_i,\phi^{NN}(t_i))\theta^{NN}(t_i)
    \right|^2,
    \\
    \operatorname{Loss}_{\mathrm{OM}}
    &=
    \frac{1}{N_T}
    \sum_{i=1}^{N_T}
    \ell_{\mathrm{OM}}
    \bigl(
    t_i,\phi^{NN}(t_i),\theta^{NN}(t_i)
    \bigr),
    \\
    \operatorname{Loss}_b
    &=
    \left|
    \phi^{NN}(0)-z_0
    \right|^2
    +
    \left|
    \phi^{NN}(T)-z_T
    \right|^2.
\end{align*}
The normalized sum in \(\operatorname{Loss}_{\mathrm{OM}}\) provides a
discrete approximation of the Onsager--Machlup running cost along the
approximated path, with the constant factor \(T\) absorbed into
\(\tau_{\mathrm{OM}}\). The positive constants \(\tau_{\phi}\),
\(\tau_{\mathrm{OM}}\), and \(\tau_b\) are penalty weights for the dynamical
constraints, the Onsager--Machlup running cost, and the boundary conditions,
respectively.

\subsection{Proofs of Theorems~\ref{theorem 3.1} and \ref{theorem3.2}}
\label{sec3.4}

After presenting the general theory in Sections \ref{sec3.1} and \ref{sec3.2}, we prove Theorems~\ref{theorem 3.1} and~\ref{theorem3.2} in four steps: Steps~1--3
establish Theorem~\ref{theorem 3.1}, and Step~4 proves Theorem~\ref{theorem3.2}.

\medskip
\noindent {\bf Step 1: Apply the Girsanov transformation to absorb the drift term.}\par
\smallskip

Consider the following auxiliary system under \(\mathbb{P}\):
\begin{align}\label{ASDE}
(\mathbb{P})
\begin{cases}
\mathrm{d}\tilde{X}_t
=
g(t,\tilde{X}_t,\tilde{Y}_t)\,\mathrm{d}t,
\qquad \tilde{X}_0=x_0\in\mathbb{R}^{d},\\[2mm]
\mathrm{d}\tilde{Y}_t
=
\dot{\phi}_2(t)\,\mathrm{d}t+\sigma(t)\,\mathrm{d}W_t,
\qquad \tilde{Y}_0=y_0\in\mathbb{R}^{m}.
\end{cases}
\end{align}
Define
\[
    \eta_t
    :=
    \sigma^{-1}(t)
    \bigl[
        f(t,\tilde X_t,\tilde Y_t)-\dot\phi_2(t)
    \bigr],
    \qquad t\in[0,T],
\]
and
\[
    \tilde R_t
    :=
    \exp\left\{
        \int_0^t\langle \eta_s,\mathrm{d}W_s\rangle
        -
        \frac12\int_0^t|\eta_s|^2\,\mathrm{d}s
    \right\}.
\]
By the assumptions of the theorem, \(f\), \(\sigma^{-1}\), and \(\dot\phi_2\) are bounded. Hence \(\eta\) is bounded, and Novikov's condition
is satisfied, namely 
\(
    \mathbb E_{\mathbb P}
    \left[
    \exp\left\{
        \frac12\int_0^T |\eta_t|^2\,\mathrm{d}t
    \right\}
    \right]
    <\infty .
\)
Therefore \((\tilde R_t)_{t\in[0,T]}\) is a martingale and \(\mathbb E_{\mathbb P}\tilde R_T=1\). Thus
\[
    \mathbb Q(A)
    :=
    \int_A \tilde R_T(\omega)\,\mathrm d\mathbb P(\omega),
    \qquad A\in\mathcal F_T,
\]
defines a probability measure on \((\Omega,\mathcal F_T)\). Equivalently,
\(
    \frac{\mathrm d\mathbb Q}{\mathrm d\mathbb P}
    =
    \tilde R_T .
\)
By Girsanov's theorem,
\[
    \tilde W_t
    :=
    W_t-\int_0^t \eta_s\,\mathrm d s
\]
is an \(m\)-dimensional Brownian motion under \(\mathbb Q\). Hence, under \(\mathbb Q\), we have
\(
    \mathrm dW_t
    =
    \mathrm d\tilde W_t+\eta_t\,\mathrm d t .
\)
Substituting this into \eqref{ASDE}, we obtain
\[
(\mathbb Q)
\begin{cases}
\mathrm{d}\tilde X_t
=
g(t,\tilde X_t,\tilde Y_t)\,\mathrm{d}t,
\qquad \tilde X_0=x_0,\\[2mm]
\mathrm{d}\tilde Y_t
=
f(t,\tilde X_t,\tilde Y_t)\,\mathrm{d}t
+
\sigma(t)\,\mathrm{d}\tilde W_t,
\qquad \tilde Y_0=y_0 .
\end{cases}
\]
Consequently, by uniqueness in law,
\(
    \mathcal L_{\mathbb Q}(\tilde X,\tilde Y)
    =
    \mathcal L_{\mathbb P}(X,Y).
\)

Let
\[
    Z_t=
    \begin{pmatrix}
        X_t\\
        Y_t
    \end{pmatrix},
    \quad
    \tilde Z_t=
    \begin{pmatrix}
        \tilde X_t\\
        \tilde Y_t
    \end{pmatrix},
    \quad
    \phi(t)=
    \begin{pmatrix}
        \phi_1(t)\\
        \phi_2(t)
    \end{pmatrix}.
\]
By the change of measure obtained above, we
have
\begin{align}\label{conclusion0}
    \mathbb P\bigl(\|Z-\phi\|<\varepsilon\bigr)
    &=
    \mathbb Q\bigl(\|\tilde Z-\phi\|<\varepsilon\bigr) \notag\\
    &=
    \mathbb E_{\mathbb P}
    \left[
        \tilde R_T
        \mathbf 1_{\{\|\tilde Z-\phi\|<\varepsilon\}}
    \right] \notag\\
    &=
    \mathbb E_{\mathbb P}
    \left[
        \tilde R_T
        \,\middle|\,
        \|\tilde Z-\phi\|<\varepsilon
    \right]
    \mathbb P\bigl(\|\tilde Z-\phi\|<\varepsilon\bigr).
\end{align}
The invertibility and regularity of \(\sigma(t)\) 
imply that \(\tilde Z-\phi\) and \(W\) are equivalent in the corresponding path norm. More precisely, there exist constants \(c_1,c_2>0\), independent of \(\varepsilon\), such that
\begin{equation}\label{conclusion1}
    c_1\|W\|
    \leq
    \|\tilde Z-\phi\|
    \leq
    c_2\|W\|.
\end{equation}
Equivalently,
\(
    \|\tilde Z-\phi\|=O(\varepsilon)
    \Longleftrightarrow
    \|W\|=O(\varepsilon).
\)
A detailed proof is given in Appendix.

It follows from \eqref{conclusion1} that
\[
    \left\{
        \|W\|<\frac{\varepsilon}{c_2}
    \right\}
    \subset
    \left\{
        \|\tilde Z-\phi\|<\varepsilon
    \right\}
    \subset
    \left\{
        \|W\|<\frac{\varepsilon}{c_1}
    \right\}.
\]
Hence, up to a multiplicative change of the radius of the small ball, the conditioning event \(\{\|\tilde Z-\phi\|<\varepsilon\}\) is equivalent to the Brownian small-ball event \(\{\|W\|<\varepsilon\}\). Therefore, for the purpose of the small-ball asymptotic analysis, \eqref{conclusion0} may be written as
\begin{align}\label{conclusion2}
   \frac{
        \mathbb P\bigl(\|Z-\phi\|<\varepsilon\bigr)}
        {\mathbb P\bigl(\|W\|<\varepsilon\bigr)}
   &\asymp
   \mathbb E_{\mathbb P}
   \left[
        \exp\left\{
            \int_0^T
            \langle \eta_t,\mathrm dW_t\rangle
            -
            \frac12
            \int_0^T
            |\eta_t|^2\,\mathrm dt
        \right\}
        \,\middle|\,
        \|W\|<\varepsilon
   \right],
   \qquad \varepsilon\downarrow0 .
\end{align}

\noindent {\bf Step 2: Reorganize the integral terms}\par
\smallskip

The exponent in \eqref{conclusion2} consists of two integral terms. We first consider the stochastic integral
\begin{equation*}
   \int_0^T \langle \eta_t,\mathrm dW_t\rangle
  =
   \int_0^T
   \left\langle
        \sigma^{-1}(t)
        \bigl[
            f(t,\tilde X_t,\tilde Y_t)-\dot\phi_2(t)
        \bigr],
        \mathrm dW_t
   \right\rangle  
   \triangleq I_1+I_2,
\end{equation*}
where
\begin{align*}
    I_1
    &:=
    \int_0^T
    \left\langle
        \sigma^{-1}(t)f(t,\tilde X_t,\tilde Y_t),
        \mathrm dW_t
    \right\rangle,\\
    I_2
    &:=
    -
    \int_0^T
    \left\langle
        \sigma^{-1}(t)\dot\phi_2(t),
        \mathrm dW_t
    \right\rangle .
\end{align*}
Write
\(
    \sigma^{-1}(t)=\bigl(a_{ij}(t)\bigr)_{m\times m},
    f=(f_1,\ldots,f_m)^T .
\)
Then
\begin{align*}
    I_1
    =
    \sum_{i=1}^m
    \int_0^T
    \sum_{j=1}^m
    a_{ij}(t)
    f_j(t,\tilde X_t,\tilde Y_t)
    \,\mathrm dW_t^i .
\end{align*}
Here \(W_t^i\) denotes the \(i\)-th component of the \(m\)-dimensional standard Brownian motion. Define the deviation processes
\[
    \hat W_t
    =
    \bigl(\hat W_t^1,\ldots,\hat W_t^d\bigr)^T,
    \quad
    \overline W_t
    =
    \bigl(\overline W_t^1,\ldots,\overline W_t^m\bigr)^T,
\]
where
\(
    \hat W_t^l
    :=
    \tilde X_t^l-\phi_1^l(t),
    \overline W_t^k
    :=
    \tilde Y_t^k-\phi_2^k(t), \mathrm d\overline W_t^k
    =
    \sum_{q=1}^m
    \sigma_{kq}(t)\,\mathrm dW_t^q, k=1,\ldots,m.\)
Expanding \(f_j(t,\tilde X_t,\tilde Y_t)\) around \((\phi_1(t),\phi_2(t))\), we obtain
\begin{align*}
    f_j(t,\tilde X_t,\tilde Y_t)
    &=
    f_j(t,\phi(t))
    +
    \sum_{l=1}^d
    \frac{\partial f_j}{\partial x_l}(t,\phi(t))\hat W_t^l
    +
    \sum_{k=1}^m
    \frac{\partial f_j}{\partial y_k}(t,\phi(t))\overline W_t^k
    +
    \Psi_j(t),
\end{align*}
where \(\Psi_j(t)\) denotes the higher-order remainder term. Consequently, \(I_1\) can be decomposed as
\begin{align*}
    I_1
    &=
    \sum_{i=1}^m
    \int_0^T
    \sum_{j=1}^m
    a_{ij}(t)f_j(t,\phi(t))
    \,\mathrm dW_t^i  
    +
    \sum_{i=1}^m
    \int_0^T
    \sum_{j=1}^m
    \sum_{l=1}^d
    a_{ij}(t)
    \frac{\partial f_j}{\partial x_l}(t,\phi(t))
    \hat W_t^l
    \,\mathrm dW_t^i  \\
    &\quad
    +
    \sum_{i=1}^m
    \int_0^T
    \sum_{j=1}^m
    \sum_{k=1}^m
    a_{ij}(t)
    \frac{\partial f_j}{\partial y_k}(t,\phi(t))
    \overline W_t^k
    \,\mathrm dW_t^i 
    +
    \sum_{i=1}^m
    \int_0^T
    \sum_{j=1}^m
    a_{ij}(t)\Psi_j(t)
    \,\mathrm dW_t^i \\
    &\triangleq
    I_{11}+I_{12}+I_{13}+I_{14}.
\end{align*}

Next, we reorganize the quadratic term in \eqref{conclusion2}. Let
\(
    R(t):=(\sigma^{-1}(t))^T\sigma^{-1}(t)
    =
    \bigl(r_{ij}(t)\bigr)_{m\times m}.
\)
Then
\begin{align*}
    &-\frac12
    \int_0^T
    \left|
        \sigma^{-1}(t)
        \bigl[
            f(t,\tilde X_t,\tilde Y_t)-\dot\phi_2(t)
        \bigr]
    \right|^2
    \,\mathrm dt  \\
    &=
    -\frac12
    \int_0^T
    \sum_{i=1}^m
    \sum_{j=1}^m
    r_{ij}(t)
    \bigl[
        f_i(t,\tilde X_t,\tilde Y_t)-\dot\phi_2^i(t)
    \bigr]
    \bigl[
        f_j(t,\tilde X_t,\tilde Y_t)-\dot\phi_2^j(t)
    \bigr]
    \,\mathrm dt \\
    &\triangleq
    -\frac12
    \int_0^T
    \bigl(M_1(t)+M_2(t)\bigr)\,\mathrm dt ,
\end{align*}
where the diagonal and off-diagonal parts are defined by
\begin{align*}
    M_1(t)
    &:=
    \sum_{i=1}^m
    r_{ii}(t)
    \bigl[
        f_i(t,\tilde X_t,\tilde Y_t)-\dot\phi_2^i(t)
    \bigr]^2,\\
    M_2(t)
    &:=
    \sum_{i\neq j}
    r_{ij}(t)
    \bigl[
        f_i(t,\tilde X_t,\tilde Y_t)-\dot\phi_2^i(t)
    \bigr]
    \bigl[
        f_j(t,\tilde X_t,\tilde Y_t)-\dot\phi_2^j(t)
    \bigr].
\end{align*}
Here \(\phi_2^i(t)\) denotes the \(i\)-th component of \(\phi_2(t)\).

For \(M_1(t)\), we write
\begin{align*}
    M_1(t)
    &=
    \sum_{i=1}^m
    r_{ii}(t)
    \bigl[
        f_i(t,\tilde X_t,\tilde Y_t)
        -
        f_i(t,\phi(t))
        +
        f_i(t,\phi(t))
        -
        \dot\phi_2^i(t)
    \bigr]^2  \\
    &=
    \sum_{i=1}^m
    r_{ii}(t)
    \bigl[
        f_i(t,\phi(t))-\dot\phi_2^i(t)
    \bigr]^2 
    +
    \sum_{i=1}^m
    r_{ii}(t)
    \bigl[
        f_i^2(t,\tilde X_t,\tilde Y_t)
        -
        f_i^2(t,\phi(t))
    \bigr]  \\
    &\quad
    +
    2
    \sum_{i=1}^m
    r_{ii}(t)
    \bigl[
        f_i(t,\phi(t))
        -
        f_i(t,\tilde X_t,\tilde Y_t)
    \bigr]
    \dot\phi_2^i(t) \\
    &\triangleq
    M_{11}(t)+M_{12}(t)+M_{13}(t).
\end{align*}

Similarly, the off-diagonal part can be decomposed as
\begin{align*}
    M_2(t)
    &=
    \sum_{i\neq j}
    r_{ij}(t)
    \Bigl[
        f_i(t,\tilde X_t,\tilde Y_t)-f_i(t,\phi(t))
        +
        f_i(t,\phi(t))-\dot\phi_2^i(t)
    \Bigr] \\
    &\quad \times
    \Bigl[
        f_j(t,\tilde X_t,\tilde Y_t)-f_j(t,\phi(t))
        +
        f_j(t,\phi(t))-\dot\phi_2^j(t)
    \Bigr] \\
    &=
    \sum_{i\neq j}
    r_{ij}(t)
    \bigl[
        f_i(t,\phi(t))-\dot\phi_2^i(t)
    \bigr]
    \bigl[
        f_j(t,\phi(t))-\dot\phi_2^j(t)
    \bigr] \\
    &\quad
    +
    \sum_{i\neq j}
    r_{ij}(t)
    \bigl[
        f_i(t,\tilde X_t,\tilde Y_t)-f_i(t,\phi(t))
    \bigr]
    \bigl[
        f_j(t,\tilde X_t,\tilde Y_t)-f_j(t,\phi(t))
    \bigr] \\
    &\quad
    +
    \sum_{i\neq j}
    r_{ij}(t)
    \Bigl\{
    \bigl[
        f_i(t,\tilde X_t,\tilde Y_t)-f_i(t,\phi(t))
    \bigr]
    \bigl[
        f_j(t,\phi(t))-\dot\phi_2^j(t)
    \bigr] \\
    &\qquad\qquad
    +
    \bigl[
        f_j(t,\tilde X_t,\tilde Y_t)-f_j(t,\phi(t))
    \bigr]
    \bigl[
        f_i(t,\phi(t))-\dot\phi_2^i(t)
    \bigr]
    \Bigr\} \\
    &\triangleq
    M_{21}(t)+M_{22}(t)+M_{23}(t).
\end{align*}

\noindent {\bf Step 3: Extract the Onsager--Machlup action functional from the individual terms}\par
\smallskip

We now identify the contribution of each term in the exponential. 
For \(I_{11}\), Lemma~\ref{lem:borell_ledoux} gives
\begin{align}\label{conclusionI11}
    \lim_{\varepsilon\to 0}
    \mathbb{E}\left(
        \exp(I_{11})\,\middle|\, \|W\|<\varepsilon
    \right)
    =1 .
\end{align}

Next, for \(I_{12}\), set
\[
    b_{ijl}(t)
    :=
    a_{ij}(t)
    \frac{\partial f_j}{\partial x_l}(t,\phi(t)).
\]
Then
\[
    I_{12}
    =
    \sum_{i=1}^m\sum_{j=1}^m\sum_{l=1}^d
    \int_0^T b_{ijl}(t)\hat W_t^l\,\mathrm dW_t^i .
\]
Applying It\^{o}'s formula to
\(b_{ijl}(t)\hat W_t^l W_t^i\), we obtain
\begin{align*}
    \mathrm d\bigl[
        b_{ijl}(t)\hat W_t^l W_t^i
    \bigr]
    &=
    b'_{ijl}(t)\hat W_t^l W_t^i\,\mathrm dt
    +
    b_{ijl}(t)\hat W_t^l\,\mathrm dW_t^i 
    +
    b_{ijl}(t)W_t^i\,\mathrm d\hat W_t^l
    +
    b_{ijl}(t)\,\mathrm d\langle \hat W^l,W^i\rangle_t .
\end{align*}
Since
\[
    \mathrm d\hat W_t^l
    =
    \bigl[
        g_l(t,\tilde X_t,\tilde Y_t)-\dot\phi_1^l(t)
    \bigr]\,\mathrm dt,
    \quad
    \mathrm d\langle \hat W^l,W^i\rangle_t=0,
\]
it follows that
\begin{align*}
    I_{12}
    &=
    \sum_{i=1}^m\sum_{j=1}^m\sum_{l=1}^d
    \Bigg\{
        b_{ijl}(T)\hat W_T^l W_T^i   
        -
        \int_0^T b'_{ijl}(t)\hat W_t^l W_t^i\,\mathrm dt
       - \int_0^T b_{ijl}(t)W_t^i
        \bigl[
            g_l(t,\tilde X_t,\tilde Y_t)-\dot\phi_1^l(t)
        \bigr]\,\mathrm dt
    \Bigg\}.
\end{align*}
On the event \(\{\|W\|<\varepsilon\}\), both \(\hat W_t\) and \(W_t\) are uniformly small. Hence, by the Lipschitz continuity of \(g\) and the
boundedness of \(\phi\), \(\dot\phi\), and \(b_{ijl}\), the preceding identity
implies
\begin{align}\label{conclusionI12}
    \lim_{\varepsilon\to 0}
    \mathbb{E}\left(
        \exp(I_{12})\,\middle|\, \|W\|<\varepsilon
    \right)
    =1 .
\end{align}

For \(I_{13}\), define
\[
    \alpha_{ijk}(t)
    :=
    a_{ij}(t)
    \frac{\partial f_j}{\partial y_k}(t,\phi(t)).
\]
Then
\[
    I_{13}
    =
    \sum_{i=1}^m\sum_{j=1}^m\sum_{k=1}^m
    \int_0^T \alpha_{ijk}(t)\overline W_t^k\,\mathrm dW_t^i .
\]
Applying It\^{o}'s formula to
\(\alpha_{ijk}(t)\overline W_t^k W_t^i\), we have
\begin{align*}
    \mathrm d\bigl[
        \alpha_{ijk}(t)\overline W_t^k W_t^i
    \bigr]
    &=
    \alpha'_{ijk}(t)\overline W_t^k W_t^i\,\mathrm dt
    +
    \alpha_{ijk}(t)\overline W_t^k\,\mathrm dW_t^i  \\
    &\quad
    +
    \alpha_{ijk}(t)W_t^i\,\mathrm d\overline W_t^k
    +
    \alpha_{ijk}(t)\,\mathrm d\langle \overline W^k,W^i\rangle_t .
\end{align*}
Since
\[
    \mathrm d\overline W_t^k
    =
    \sum_{q=1}^m\sigma_{kq}(t)\,\mathrm dW_t^q,
    \quad
    \mathrm d\langle \overline W^k,W^i\rangle_t
    =
    \sigma_{ki}(t)\,\mathrm dt,
\]
we obtain
\begin{align*}
    I_{13}
    &=
    \sum_{i=1}^m\sum_{j=1}^m\sum_{k=1}^m
    \Bigg\{
        \alpha_{ijk}(T)\overline W_T^k W_T^i
        -
        \int_0^T
        \alpha'_{ijk}(t)\overline W_t^k W_t^i\,\mathrm dt
        -
        \int_0^T
        \alpha_{ijk}(t)W_t^i\,\mathrm d\overline W_t^k  \\
    &\qquad
        -
        \int_0^T
        \alpha_{ijk}(t)\sigma_{ki}(t)\,\mathrm dt
    \Bigg\}  \\
    &\triangleq
    \sum_{i=1}^m\sum_{j=1}^m\sum_{k=1}^m
    \bigl(
        L_{ijk}^1
        -
        L_{ijk}^2
        -
        L_{ijk}^3
        -
        L_{ijk}^4
    \bigr).
\end{align*}
On the event \(\{\|W\|<\varepsilon\}\), both \(W_t\) and \(\overline W_t\) are uniformly small. Therefore,
\begin{align}\label{conclusionL1L2}
    \lim_{\varepsilon\to 0}
    \mathbb{E}\left(
        \exp(L_{ijk}^1)\,\middle|\, \|W\|<\varepsilon
    \right)
    =
    \lim_{\varepsilon\to 0}
    \mathbb{E}\left(
        \exp(L_{ijk}^2)\,\middle|\, \|W\|<\varepsilon
    \right)
    =1 .
\end{align}
For \(L_{ijk}^3\), we have
\begin{equation*}
    L_{ijk}^3=
    \int_0^T
    \alpha_{ijk}(t)W_t^i\,\mathrm d\overline W_t^k=
    \sum_{q=1}^m
    \int_0^T
    \alpha_{ijk}(t)\sigma_{kq}(t)W_t^i\,\mathrm dW_t^q .
\end{equation*}
For the terms with \(q\neq i\), Lemma~\ref{lem:shepp_zeitouni} gives
\begin{align*}
    \lim_{\epsilon\to0}
    \mathbb{E}\left[
        \exp\left(
            \int_0^T
            \alpha_{ijk}(t)\sigma_{kq}(t)W_t^i\,\mathrm dW_t^q
        \right)
        \,\middle|\, \|W\|<\varepsilon
    \right]
    =
    1 .
\end{align*}
For the term with \(q=i\), It\^{o}'s formula gives
\begin{align*}
    \int_0^T
    \alpha_{ijk}(t)\sigma_{ki}(t)W_t^i\,\mathrm dW_t^i
    =&
    \frac12\alpha_{ijk}(T)\sigma_{ki}(T)(W_T^i)^2
    -
    \frac12\int_0^T
    \bigl(\alpha_{ijk}(t)\sigma_{ki}(t)\bigr)'(W_t^i)^2\,\mathrm dt
   \\ &-
    \frac12\int_0^T
    \alpha_{ijk}(t)\sigma_{ki}(t)\,\mathrm dt .
\end{align*}
Hence, again using the small-ball conditioning, we obtain
\begin{align}\label{conclusionL3}
    \lim_{\varepsilon\to0}
    \mathbb{E}\left[
        \exp\left(
            L_{ijk}^3
            +
            \frac12
            \int_0^T
            \alpha_{ijk}(t)\sigma_{ki}(t)\,\mathrm dt
        \right)
        \,\middle|\, \|W\|<\varepsilon
    \right]
    =
    1 .
\end{align}

For the remainder term \(I_{14}\), following the argument in \cite{liu_onsagermachlup_2024}, one can prove that
\begin{align}\label{conclusionI14}
    \lim_{\varepsilon\to 0}
    \mathbb{E}\left(
        \exp(I_{14})\,\middle|\, \|W\|<\varepsilon
    \right)
    =1 .
\end{align}

Combining \eqref{conclusionI11}, \eqref{conclusionI12}, \eqref{conclusionL1L2}, \eqref{conclusionL3}, and \eqref{conclusionI14}, we obtain
\begin{align*}
    \lim_{\varepsilon\to0}
    \mathbb{E}\left[
        \exp\left(
            I_1
            +
            \frac12
            \sum_{i=1}^m\sum_{j=1}^m\sum_{k=1}^m
            \int_0^T
            a_{ij}(t)
            \frac{\partial f_j}{\partial y_k}(t,\phi(t))
            \sigma_{ki}(t)\,\mathrm dt
        \right)
        \,\middle|\, \|W\|<\varepsilon
    \right]
    =
    1 .
\end{align*}
Since \(a(t)=\sigma^{-1}(t)\), we have
\[
    \sum_{i=1}^m a_{ij}(t)\sigma_{ki}(t)
    =
    \delta_{jk}.
\]
Therefore,
\[
    \sum_{i=1}^m\sum_{j=1}^m\sum_{k=1}^m
    a_{ij}(t)
    \frac{\partial f_j}{\partial y_k}(t,\phi(t))
    \sigma_{ki}(t)
    =
    \sum_{k=1}^m
    \frac{\partial f_k}{\partial y_k}(t,\phi(t))
    =
    \operatorname{div}_y f(t,\phi(t)).
\]
Thus,
\begin{align}\label{conclusionI1}
    \lim_{\varepsilon\to0}
    \mathbb{E}\left[
        \exp\left(
            I_1
            +
            \frac12
            \int_0^T
            \operatorname{div}_y f(t,\phi(t))\,\mathrm dt
        \right)
        \,\middle|\, \|W\|<\varepsilon
    \right]
    =
    1 .
\end{align}

For \(I_2\), Lemma~\ref{lem:borell_ledoux} yields
\begin{align}\label{conclusionI2}
    \lim_{\varepsilon\to 0}
    \mathbb{E}\left(
        \exp(I_2)\,\middle|\, \|W\|<\varepsilon
    \right)
    =1 .
\end{align}

It remains to consider the quadratic term. By the Lipschitz continuity of \(f\), the boundedness of \(\phi\) and \(\dot\phi\), and the small-ball conditioning, the error terms in the decomposition of \(M_1\) and \(M_2\) vanish in the exponential sense. More precisely,
\begin{align*}
    \lim_{\varepsilon\to 0}
    \mathbb{E}\left[
        \exp\left(
            -\frac12
            \int_0^T
            \bigl(
                M_{12}(t)+M_{13}(t)+M_{22}(t)+M_{23}(t)
            \bigr)\,\mathrm dt
        \right)
        \,\middle|\, \|W\|<\varepsilon
    \right]
    =
    1 .
\end{align*}
The remaining leading terms satisfy
\begin{align*}
    M_{11}(t)+M_{21}(t)
    &=
    \bigl(
        f_1(t,\phi(t))-\dot\phi_2^1(t),
        \ldots,
        f_m(t,\phi(t))-\dot\phi_2^m(t)
    \bigr)  \\
    &\quad \times
    R(t)
    \bigl(
        f_1(t,\phi(t))-\dot\phi_2^1(t),
        \ldots,
        f_m(t,\phi(t))-\dot\phi_2^m(t)
    \bigr)^T  \\
    &=
    \left|
        \sigma^{-1}(t)
        \bigl[
            f(t,\phi(t))-\dot\phi_2(t)
        \bigr]
    \right|^2 .
\end{align*}
Consequently,
\begin{align}\label{conclusionM1M2}
    &\lim_{\varepsilon \to 0}
    \mathbb{E}\left[
        \exp\left(
            -\frac12
            \int_0^T
            \left|
                \sigma^{-1}(t)
                \bigl[
                    f(t,\tilde X_t,\tilde Y_t)-\dot\phi_2(t)
                \bigr]
            \right|^2
            \,\mathrm dt
        \right)
        \,\middle|\, \|W\|<\varepsilon
    \right]
    \notag\\
    &\qquad =
    \exp\left(
        -\frac12
        \int_0^T
        \left|
            \sigma^{-1}(t)
            \bigl[
                f(t,\phi(t))-\dot\phi_2(t)
            \bigr]
        \right|^2
        \,\mathrm dt
    \right).
\end{align}

Combining \eqref{conclusionI1}, \eqref{conclusionI2}, and \eqref{conclusionM1M2}, we finally obtain the small-ball asymptotic relation
\begin{align*}
    \frac{
        \mathbb{P}\bigl(\|Z-\phi\|<\varepsilon\bigr)
    }{
        \mathbb{P}\bigl(\|W\|<\varepsilon\bigr)
    }
    \sim
    \exp\Bigg\{
        -\frac12
        \int_0^T
        \left|
            \sigma^{-1}(t)
            \bigl[
                f(t,\phi(t))-\dot\phi_2(t)
            \bigr]
        \right|^2
        \,\mathrm dt
        -
        \frac12
        \int_0^T
        \operatorname{div}_y f(t,\phi(t))\,\mathrm dt
    \Bigg\}
\end{align*}
as \(\varepsilon\to0\). This completes the proof of Theorem~\ref{theorem 3.1}.
\hfill \(\qedsymbol\)

\noindent {\bf Step 4: Proof of Theorem~\ref{theorem3.2} by Theorem~\ref{theorem 3.1}.}\par
\smallskip

We now prove Theorem~\ref{theorem3.2}. Let \(Z_t=(X_t,Y_t)\). Since the
Lamperti transformation is defined on the domain \(\mathcal D\), we restrict
the reference path to satisfy
\(
    \phi([0,T])\Subset\mathcal D .
\)
Then
\(
    \delta_\phi
    :=
\operatorname{dist}\bigl(\phi([0,T]),\mathcal D^c\bigr)>0.
\)
Hence, for the uniform norm,
every \(0<\varepsilon<\delta_\phi\) gives
\(
    \{\|Z-\phi\|<\varepsilon\}
    \subset
    \{Z_t\in\mathcal D,\ 0\leq t\leq T\}.
\)
Therefore, the Lamperti transformation is valid on the small-ball events used
in the Onsager--Machlup asymptotics.

For the stochastic differential equation~\eqref{E-1}, define
\begin{align*}
    U_t
    =
    \Gamma(t,X_t,Y_t)
    =
    \bigl(\Gamma_1(t,X_t,Y_t),\ldots,\Gamma_m(t,X_t,Y_t)\bigr)^T .
\end{align*}
By It\^{o}'s formula, we have
\begin{align*}
    \mathrm dU_t
    &=
    \Bigg[
        \partial_t \Gamma(t,X_t,Y_t)
        +
        \partial_x \Gamma(t,X_t,Y_t)g(t,X_t,Y_t)
        +
        \partial_y \Gamma(t,X_t,Y_t)f(t,X_t,Y_t)
\\
    &\qquad\qquad
        +
        \frac12
        \operatorname{Tr}
        \left(
            \partial_{yy}^2 \Gamma(t,X_t,Y_t)
            \sigma\sigma^T(t,X_t,Y_t)
        \right)
    \Bigg]\,\mathrm dt
    +
    \partial_y \Gamma(t,X_t,Y_t)
    \sigma(t,X_t,Y_t)\,\mathrm dW_t .
\end{align*}
Here
\(
    \partial_t \Gamma
    :=
    \left(
        \partial_t \Gamma_1,
        \ldots,
        \partial_t \Gamma_m
    \right)^T,
    \partial_x \Gamma
    :=
    \left(
        \partial_{x_j}\Gamma_i
    \right)_{1\le i\le m,\,1\le j\le d},
    \partial_y \Gamma
    :=
    \left(
        \partial_{y_j}\Gamma_i
    \right)_{1\le i,j\le m},
\)
and
\(
    \operatorname{Tr}
    \left(
        \partial_{yy}^2 \Gamma
        \sigma\sigma^T
    \right)
    :=
    \left(
        \operatorname{tr}
        \left(
            \partial_{yy}^2 \Gamma_1
            \sigma\sigma^T
        \right),
        \ldots,
        \operatorname{tr}
        \left(
            \partial_{yy}^2 \Gamma_m
            \sigma\sigma^T
        \right)
    \right)^T .
\)
By the definition of the Lamperti transformation in
Assumption~\ref{ass:lamperti-domain},
\(
    \partial_y \Gamma(t,x,y)\sigma(t,x,y)=I_m,
\)
or equivalently,
\(
    \partial_y \Gamma(t,x,y)=\sigma^{-1}(t,x,y),
(t,x,y)\in [0,T]\times\mathcal D .
\)
Therefore,
\[
    \mathrm dU_t
    =
    H(t,X_t,Y_t)\,\mathrm dt+\mathrm dW_t,
\]
where
\begin{align*}
    H(t,x,y)
    &:=
    \partial_t \Gamma(t,x,y)
    +
    \partial_x \Gamma(t,x,y)g(t,x,y)
    +
    \sigma^{-1}(t,x,y)f(t,x,y)
\\
    &\qquad
    +
    \frac12
    \operatorname{Tr}
    \left(
        \partial_{yy}^2 \Gamma(t,x,y)
        \sigma\sigma^T(t,x,y)
    \right).
\end{align*}
By the local invertibility assumption in
Assumption~\ref{ass:lamperti-domain}, the map
\(
    y\mapsto \Gamma(t,x,y)=:u
\)
is invertible on the Lamperti domain. Its inverse with respect to the
\(y\)-variable is denoted by
\(
    y=v(t,x,u).
\)
Consequently, inside the Lamperti chart, the transformed process
\((X_t,U_t)\) satisfies
\begin{equation}\label{eq:transformed-XU-system}
    \begin{cases}
    \mathrm dX_t
    =
    \widetilde g(t,X_t,U_t)\,\mathrm dt,\\[1mm]
    \mathrm dU_t
    =
    \widetilde H(t,X_t,U_t)\,\mathrm dt+\mathrm dW_t,
    \end{cases}
\end{equation}
where
\(
    \widetilde g(t,x,u)
    :=
    g(t,x,v(t,x,u)),
    \widetilde H(t,x,u)
    :=
    H(t,x,v(t,x,u)).
\)

Define the transformed reference path by
\[
    \phi_3(t)
    :=
    u(t,\phi_1(t),\phi_2(t)).
\]
Since \(\phi([0,T])\Subset\mathcal D\), the transformed path
\((\phi_1,\phi_3)([0,T])\) lies inside the corresponding Lamperti coordinate
domain. Applying Theorem~\ref{theorem 3.1} to the transformed additive-noise
system~\eqref{eq:transformed-XU-system}, we obtain
\[
    I(\phi_1,\phi_3)
    =
    \int_0^T
    \left[
        \frac12
        \left|
            \dot{\phi}_3(t)
            -
            \widetilde H(t,\phi_1(t),\phi_3(t))
        \right|^2
        +
        \frac12
        \operatorname{div}_{u}
        \widetilde H(t,\phi_1(t),\phi_3(t))
    \right]
    \mathrm dt .
\]
Moreover, the transformed reference path satisfies the deterministic
constraint
\[
    \dot{\phi}_1(t)
    =
    \widetilde g(t,\phi_1(t),\phi_3(t))
    =
    g(t,\phi_1(t),\phi_2(t)).
\]

We now rewrite the two terms in the integrand in the original variables.
Since
\(
\phi_3(t)=\Gamma(t,\phi_1(t),\phi_2(t)),
\)
we have
\[
    \dot{\phi}_3(t)
    =
    \partial_t \Gamma(t,\phi(t))
    +
    \partial_x \Gamma(t,\phi(t))\dot{\phi}_1(t)
    +
    \partial_y \Gamma(t,\phi(t))\dot{\phi}_2(t).
\]
Using
\(
    \partial_y u(t,\phi(t))
    =
    \sigma^{-1}(t,\phi(t)),
    \dot{\phi}_1(t)=g(t,\phi(t)),
\)
we obtain
\begin{align*}
    \dot{\phi}_3(t)
    -
    \widetilde H(t,\phi_1(t),\phi_3(t))
    &=
    \sigma^{-1}(t,\phi(t))
    \bigl[
        \dot{\phi}_2(t)-f(t,\phi(t))
    \bigr]
    -
    \frac12
    \operatorname{Tr}
    \left(
        \partial_{yy}^2\Gamma(t,\phi(t))
        \sigma\sigma^T(t,\phi(t))
    \right)
\\
    &=
    \sigma^{-1}(t,\phi(t))
    \left[
        \dot{\phi}_2(t)
        -
        f(t,\phi(t))
        -
        \frac12
        \sigma(t,\phi(t))
        \operatorname{Tr}
        \left(
            \partial_{yy}^2\Gamma(t,\phi(t))
            \sigma\sigma^T(t,\phi(t))
        \right)
    \right].
\end{align*}

For the divergence term, by the chain rule,
\[
    \partial_u\widetilde H(t,x,u)
    =
    \partial_y H(t,x,v(t,x,u))\,\partial_u v(t,x,u).
\]
Since \(u=\Gamma(t,x,y)\) and \(y=v(t,x,u)\), we have
\(
    \partial_y \Gamma(t,x,y)\,\partial_uv(t,x,u)=I_m.
\)
Using \(\partial_y \Gamma=\sigma^{-1}\), it follows that
\(
    \partial_u v(t,x,u)
    =
    \sigma(t,x,v(t,x,u)).
\)
Therefore, along the path \(\phi(t)\),
\begin{align*}
    \operatorname{div}_u
    \widetilde H(t,\phi_1(t),\phi_3(t))
    &=
    \operatorname{tr}
    \left[
        \partial_u
        \widetilde H(t,\phi_1(t),\phi_3(t))
    \right]
=
    \operatorname{tr}
    \left(
        \partial_y H(t,\phi(t))
        \sigma(t,\phi(t))
    \right).
\end{align*}

Combining the above identities and using the one-to-one correspondence between
\(\phi\) and \((\phi_1,\phi_3)\) in the Lamperti chart, the
Onsager--Machlup action functional for the original process \((X_t,Y_t)\) is
\(
    I(\phi)
    =
    \int_0^T
    \operatorname{OM}(t,\phi(t),\dot{\phi}(t))\,\mathrm dt,
\)
where
\begin{align*}
    \operatorname{OM}(t,\phi(t),\dot{\phi}(t))
    &=
    \frac12
    \left|
        \sigma^{-1}(t,\phi(t))
        \left[
            \dot{\phi}_2(t)
            -
            f(t,\phi(t))
            -
            \frac12
            \sigma(t,\phi(t))
            \operatorname{Tr}
            \left(
                \partial_{yy}^2\Gamma(t,\phi(t))
                \sigma\sigma^T(t,\phi(t))
            \right)
        \right]
    \right|^2
\\
    &\qquad
    +
    \frac12
    \operatorname{tr}
    \left(
        \partial_y H(t,\phi(t))
        \sigma(t,\phi(t))
    \right).
\end{align*}
This is precisely \eqref{OM2}. The proof is complete.
\hfill \(\qedsymbol\)

\section{Numerical Experiments}\label{sec4}

In this section, we consider three examples: the kinetic Langevin system and the Sisyphus-cooling system , and a stochastic damped mechanical system. 
For each example, we reformulate the computation of the MPTP as an optimal control problem and solve the resulting problem numerically by using neural networks.

\subsection{Kinetic Langevin system}

We first consider a kinetic Langevin system as a benchmark example. Although
the diffusion coefficient in this example is constant, and hence the system
belongs to the degenerate additive-noise case, it admits an analytical MPTP.
This allows us to validate the proposed neural-network optimal control method
against a known solution before applying it to genuinely multiplicative-noise
systems.

The kinetic Langevin system is a second-order stochastic dynamical system of
the form
\begin{align*}
    \ddot{X}+\gamma\dot{X}
    =
    -V'(X)+\sqrt{\mu}\gamma^{\frac{1}{2}}\dot{W}(t).
\end{align*}
Letting \(Y=\dot{X}\), the system can be rewritten as the following degenerate
stochastic differential equation:
\begin{align}\label{KLS}
    \left\{
    \begin{array}{l}
    \mathrm{d}X_t=g(X_t,Y_t)\,\mathrm{d}t,
    \qquad X_0=x_0\in \mathbb{R},
    \\[2mm]
    \mathrm{d}Y_t=f(X_t,Y_t)\,\mathrm{d}t+c\,\mathrm{d}W_t,
    \qquad Y_0=y_0\in \mathbb{R},
    \end{array}
    \right.
\end{align}
where
\(
    g(x,y)=y,
    f(x,y)=-\gamma y-V'(x),
    c=\sqrt{\mu}\gamma^{\frac{1}{2}}.
\)
The corresponding OM function is given by
\begin{align*}
    \operatorname{OM}(\phi,\dot{\phi})
    =
    \frac{1}{2}
    \left(
        \frac{
            \dot{\phi}_2+\gamma\phi_2+V'(\phi_1)
        }{c}
    \right)^2
    -
    \frac{\gamma}{2},
    \qquad
    \dot{\phi}_1=\phi_2.
\end{align*}
The associated Hamilton--Pontryagin equations are
\begin{align}\label{HP of KLS}
    \left\{
    \begin{array}{l}
    -\dfrac{1}{c^2}
    \left(
        \ddot{\phi}_2+V''(\phi_1)\dot{\phi}_1+\gamma\dot{\phi}_2
    \right)
    +
    \dfrac{
        \dot{\phi}_2+\gamma\phi_2+V'(\phi_1)
    }{c^2}
    \gamma
    +
    \lambda
    =
    0,
    \\[3mm]
    \dfrac{
        \dot{\phi}_2+\gamma\phi_2+V'(\phi_1)
    }{c^2}
    V''(\phi_1)
    +
    \dot{\lambda}
    =
    0,
    \\[3mm]
    \dot{\phi}_1=\phi_2,
    \\[1mm]
    \phi(0)=z_0=(x_0,y_0)^T,
    \\[1mm]
    \phi(T)=z_T=(x_T,y_T)^T.
    \end{array}
    \right.
\end{align}
For the special case of the quadratic potential
\(V(x)=\frac{1}{2}x^2\), an analytical solution of~\eqref{HP of KLS}
was obtained in~\cite{WOS:001514271800001}. In this example, we reformulate the MPTP problem as an optimal control problem, solve it numerically using neural networks, and compare the numerical result with the analytical solution.

The variational problem for the MPTP of system~\eqref{KLS} is
\begin{equation*}
    \left\{
    \begin{aligned}
        \min_{\phi}\quad
        &
        \int_0^T
        \left[
        \frac{1}{2}
        \left(
            \frac{
                \dot{\phi}_2+\gamma\phi_2+V'(\phi_1)
            }{c}
        \right)^2
        -
        \frac{\gamma}{2}
        \right]
        \,\mathrm{d}t,
        \\
        \text{subject to}\quad
        &
        \dot{\phi}_1=\phi_2,
        \\
        &
        \phi(0)=z_0,\qquad
        \phi(T)=z_T.
    \end{aligned}
    \right.
\end{equation*}
Since the constant term \(-\gamma/2\) does not affect the minimizer, the above problem can be equivalently reformulated as the following optimal control problem:
\begin{equation*}
    \left\{
    \begin{aligned}
        \min_{\theta}\quad
        &
        \int_0^T |\theta(t)|^2\,\mathrm{d}t,
        \\
        \text{subject to}\quad
        &
        \dot{\phi}_1=\phi_2,
        \\
        &
        \dot{\phi}_2
        =
        -\gamma\phi_2-V'(\phi_1)
        +
        \sqrt{\mu}\gamma^{\frac{1}{2}}\theta(t),
        \\
        &
        \phi(0)=(x_{0},y_{0})^T,
        \qquad
        \phi(T)=(x_{T},y_{T})^T.
    \end{aligned}
    \right.
\end{equation*}

We use fully connected neural networks to approximate the control \(\theta(t)\)
and the path \(\phi(t)=(\phi_1(t),\phi_2(t))^T\). Following the parameter
setting used in~\cite{WOS:001514271800001}, we take \(T=2\), \(N_T=200\),
\(z_0=(-1,5.2846)^T\), and \(z_T=(1,0.3099)^T\). The loss function used for
training is defined as
\(
    \operatorname{Loss}
    =
    \operatorname{Loss}_{\phi}
    +
    \operatorname{Loss}_{\theta}
    +
    10\operatorname{Loss}_{b},
\)
where
\begin{align*}
    \operatorname{Loss}_{\phi}
    =
    \frac{1}{N_T}
    \sum_{i=1}^{N_T}
    \Bigg[
        &
        \left(
            \dot{\phi}_1^{NN}(t_i)
            -
            \phi_2^{NN}(t_i)
        \right)^2
        \\
        &+
        \left(
            \dot{\phi}_2^{NN}(t_i)
            +
            \gamma \phi_2^{NN}(t_i)
            +
            V'\bigl(\phi_1^{NN}(t_i)\bigr)
            -
            \sqrt{\mu}\gamma^{\frac{1}{2}}\theta^{NN}(t_i)
        \right)^2
    \Bigg],
\end{align*}
\[\operatorname{Loss}_{\theta}
    =
    \frac{1}{N_T}
    \sum_{i=1}^{N_T}
    \left(
        \theta^{NN}(t_i)
    \right)^2,
\]
and
\begin{equation*}
    \operatorname{Loss}_{b}
    =
    \left(
        \phi_1^{NN}(0)-x_{0}
    \right)^2
    +
    \left(
        \phi_1^{NN}(T)-x_{T}
    \right)^2+
    \left(
        \phi_2^{NN}(0)-y_{0}
    \right)^2
    +
    \left(
        \phi_2^{NN}(T)-y_{T}
    \right)^2.
\end{equation*}

\begin{figure}[H]
    \centering
    \includegraphics[width=0.5\linewidth]{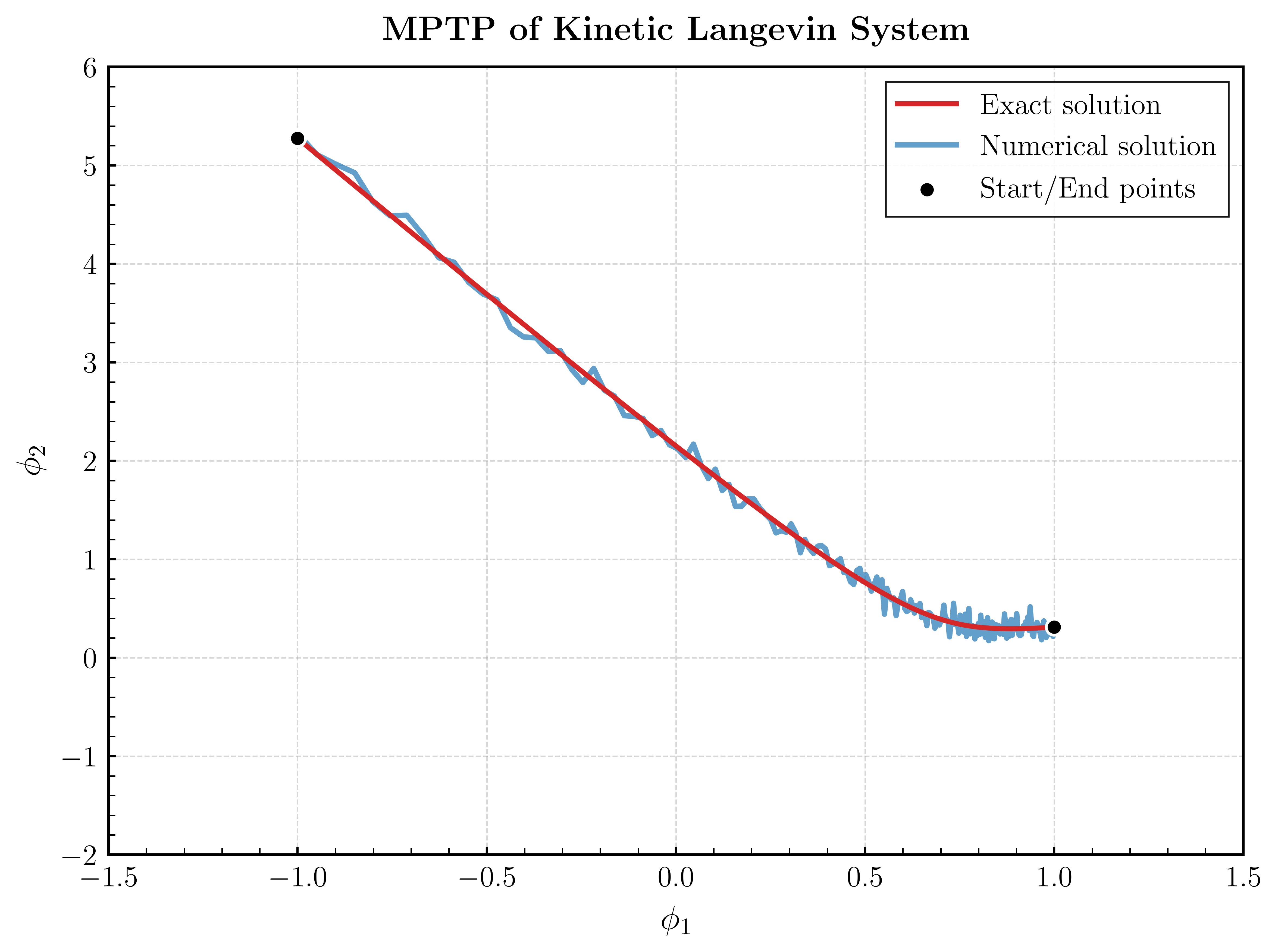}
    \caption{
    Comparison between the analytical and numerical MPTPs for  \(\gamma=3\), \(\mu=0.8\), and \(V(x)=-\frac{1}{2}x^2\). The red curve represents the analytical solution, while the blue curve represents the numerical solution obtained by the proposed method.
    }
    \label{fig:KLSfi}
\end{figure}

\begin{figure}[H]
    \centering
    \begin{subfigure}{0.47\textwidth}
        \includegraphics[width=\textwidth]{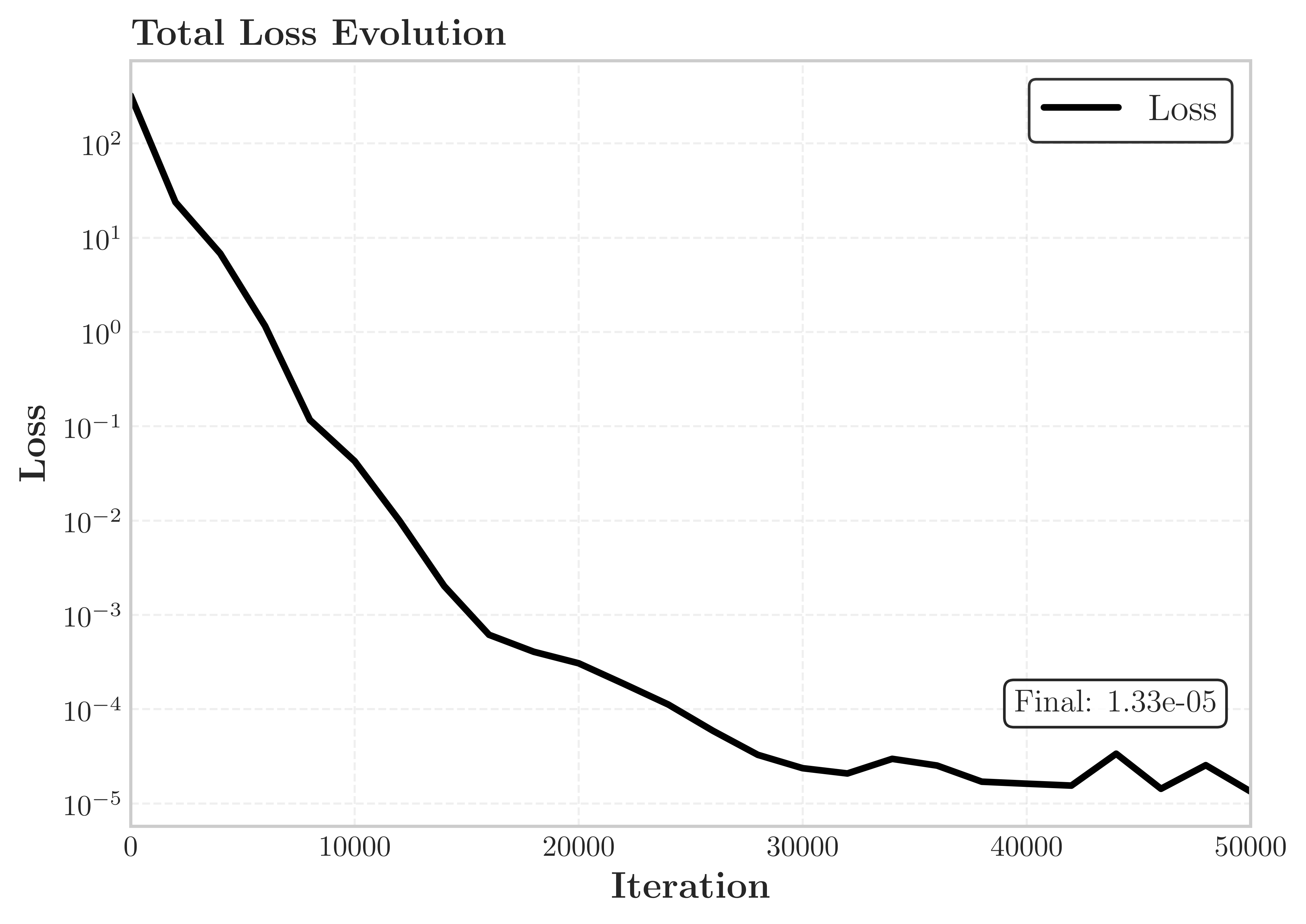}
        \caption{}
        \label{fig:1a}
    \end{subfigure}
    \begin{subfigure}{0.47\textwidth}
        \includegraphics[width=\textwidth]{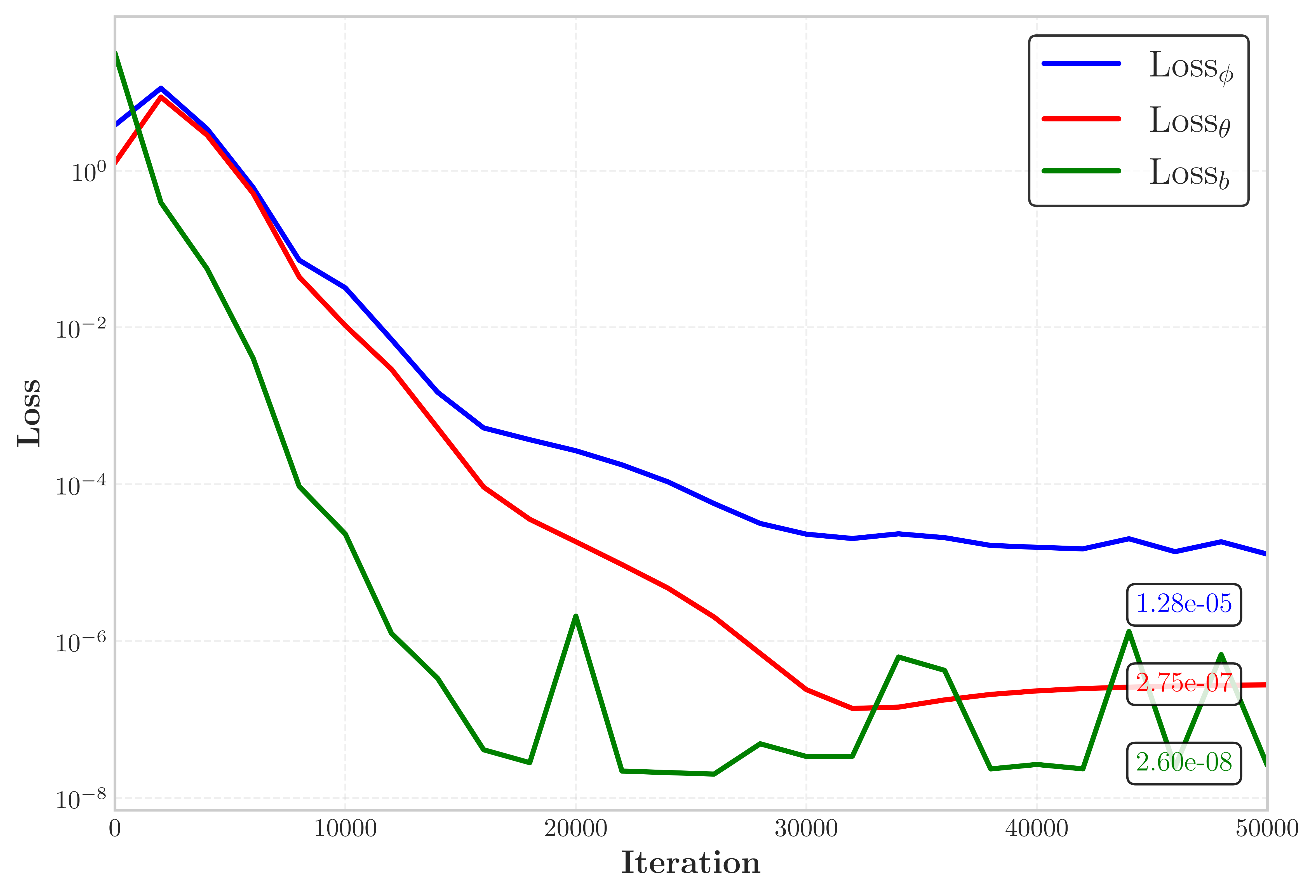}
        \caption{}
        \label{fig:1b}
    \end{subfigure}
    \caption{
    Loss curves during training.
    (a) The total loss decreases continuously during the iterations and
    eventually reaches \(1.33\times 10^{-5}\).
    (b) The blue, red, and green curves represent
    \(\operatorname{Loss}_{\phi}\), \(\operatorname{Loss}_{\theta}\), and
    \(\operatorname{Loss}_{b}\), respectively.
    }
    \label{loss}
\end{figure}

Figure~\ref{fig:KLSfi} compares the analytical MPTP and the numerical MPTP for \(\gamma=3\), \(\mu=0.8\), and \(V(x)=-\frac{1}{2}x^2\), and the transition from
\(X(0)=-1\) to \(X(2)=1\). The two curves agree well, showing that
the proposed neural-network optimal control method accurately recovers the
analytical MPTP in this benchmark additive-noise case. Figure~\ref{loss} shows the loss curves during training. The term \(\operatorname{Loss}_{\phi}\) enforces the dynamical constraints of the optimal control problem, \(\operatorname{Loss}_{\theta}\) corresponds to the objective functional, and \(\operatorname{Loss}_{b}\) penalizes the boundary conditions. At the end of the training process, all three loss terms decrease to below the order of \(10^{-5}\), which further demonstrates the numerical feasibility of the proposed method.

\subsection{The Sisyphus-cooling system driven by multiplicative Gaussian noise}

We next consider a Sisyphus-cooling system driven by multiplicative and additive
Gaussian white noises \cite{Falasco2022}. This example
is used to illustrate the application of Theorem~\ref{theorem3.2} to a
degenerate system with a genuinely state-dependent diffusion coefficient.
The system is given by
\begin{equation}
\left\{
\begin{aligned}
\mathrm{d}X_t
&=
Y_t\,\mathrm{d}t,\\
\mathrm{d}Y_t
&=
\left[
-V'(X_t)
+F(Y_t)
+D'(Y_t)
\right]\mathrm{d}t
+
\sqrt{2D(Y_t)}\,\mathrm{d}W_t,
\end{aligned}
\right.
\label{eq:sisyphus_system}
\end{equation}
where
\(
F(y)=-\frac{y}{1+y^2},
D(y)=D_0+\frac{D_1}{1+y^2},
\)
and \(W_t\) is a one-dimensional Brownian motion.

To construct a metastable transition problem, we choose the bistable potential
\(
V(x)=\frac{1}{4}x^4-\frac{1}{2}x^2.
\)
Then system~\eqref{eq:sisyphus_system} is of the form~\eqref{E-1} with
\(
g(x,y)=y,
\)
\(
f(x,y)
=
-V'(x)+F(y)+D'(y)
=
x-x^3
-\frac{y}{1+y^2}
-\frac{2D_1y}{(1+y^2)^2},
\)
and
\(
\sigma(y)
=
\sqrt{
2D(y)
}
=
\sqrt{
2\left(
D_0+\frac{D_1}{1+y^2}
\right)
}.
\)
For \(D_0>0\) and \(D_1>0\), we have
\[
\sqrt{2D_0}
\leq
\sigma(y)
\leq
\sqrt{2(D_0+D_1)},
\]
so that \(\sigma\) is uniformly nondegenerate. Define the Lamperti-type transformation
\[
\Gamma(y)
=
\int_0^y
\frac{\mathrm{d}r}{\sqrt{2D(r)}}.
\]
Then
\(
\Gamma_y(y)
=
\frac{1}{\sigma(y)}
=
\sigma^{-1}(y).
\)
Moreover, since
\(
\frac{1}{\sqrt{2(D_0+D_1)}}
\leq
\Gamma_y(y)
\leq
\frac{1}{\sqrt{2D_0}},
\)
the map \(y\mapsto \Gamma(y)\) is strictly increasing and satisfies
\(
\lim_{y\to+\infty}\Gamma(y)=+\infty,
\lim_{y\to-\infty}\Gamma(y)=-\infty.
\)
Hence the Lamperti-type transformation
required in Theorem~\ref{theorem3.2} is well defined.

According to Theorem~\ref{theorem3.2} and Remark \ref{re-scalar-case}, the OM function of
system~\eqref{eq:sisyphus_system} is
\begin{equation}
\begin{aligned}
\operatorname{OM}(\phi,\dot{\phi})
&=
\frac{
\left[
\dot{\phi}_2
+V'(\phi_1)
-F(\phi_2)
-\frac{1}{2}D'(\phi_2)
\right]^2
}{
4D(\phi_2)
}
+
\frac12
    \operatorname{tr}
    \left(
        \partial_y H(t,\phi)
        \sigma(\phi_2)
    \right),
\end{aligned}
\label{eq:sisyphus_OM}
\end{equation}
where
\(
\operatorname{tr}
    \left(
        \partial_y H(t,\phi)
        \sigma(\phi_2)
    \right)
=
F'(y)
+\frac{1}{2}D''(y)
-
\frac{
\left[
-U'(x)+F(y)+\frac{1}{2}D'(y)
\right]D'(y)
}{
2D(y)
}.
\)
Here \(\phi=(\phi_1,\phi_2)^T\) corresponds to \((X,Y)^T\).
The drift dynamics associated with \eqref{eq:sisyphus_system} has two stable
equilibria \((-1,0)\) and \((1,0)\), separated by the saddle point \((0,0)\).
We therefore consider the transition from
\(
z_0=(-1,0)
\) to \(
z_T=(1,0).
\)

By introducing the control \(\theta(t)\) through
\[
\dot{\phi}_2
=
-V'(\phi_1)
+F(\phi_2)
+D'(\phi_2)
+
\sqrt{2D(\phi_2)}\,\theta(t),
\]
the problem can be reformulated as the following optimal control problem:
\begin{equation}
\left\{
\begin{aligned}
\min_{\theta,\phi}\quad
&
\int_0^T
\left[
\frac{1}{2}
\left(
\theta(t)
+
\frac{D'(\phi_2(t))}
{2\sqrt{2D(\phi_2(t))}}
\right)^2
+
\frac12
    \operatorname{tr}
    (
        \partial_y H(t,\phi)
        \sigma(\phi_2)
    )\right],
\\[2mm]
\text{subject to}\quad
&
\dot{\phi}_1
=
\phi_2,
\\
&
\dot{\phi}_2
=
\phi_1-\phi_1^3
-\frac{\phi_2}{1+\phi_2^2}
-\frac{2D_1\phi_2}{(1+\phi_2^2)^2}
+
\sqrt{
2\left(
D_0+\frac{D_1}{1+\phi_2^2}
\right)
}
\,\theta(t),
\\
&
\phi(0)=z_0,
\qquad
\phi(T)=z_T.
\end{aligned}
\right.
\label{eq:sisyphus_optimal_control}
\end{equation}

\begin{figure}[H]
\centering
\includegraphics[width=0.95\textwidth]{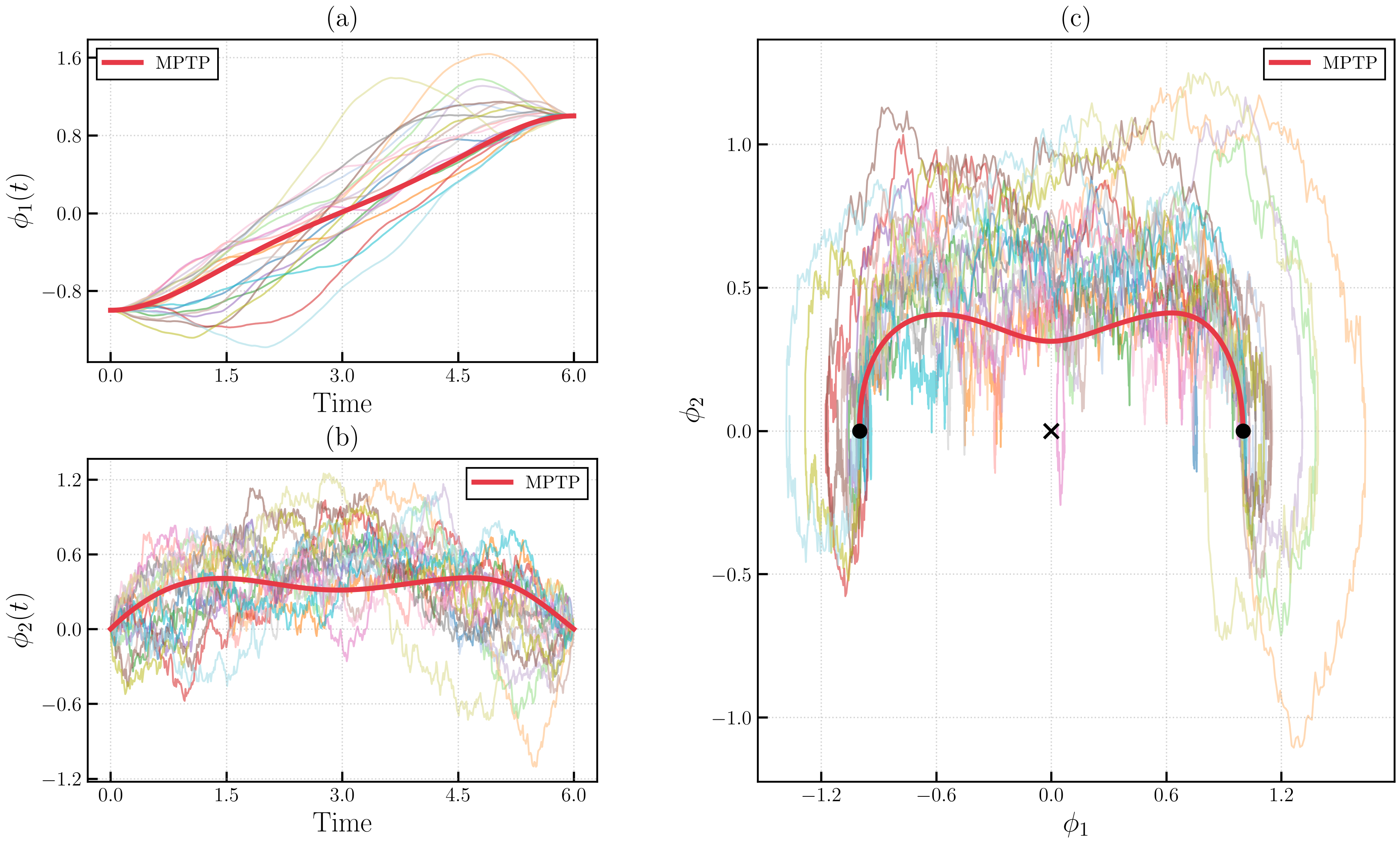}
\caption{
Numerical MPTP and simulated transition trajectories for the bistable
Sisyphus-cooling system driven by multiplicative Gaussian noise.
The red curve represents the numerical MPTP obtained by the proposed method.
}
\label{fig:sisyphus_mptp}
\end{figure}

For
\(
D_0=0.02, 
D_1=0.08,
T=6,
\)
we compute the MPTP from
\(
z_0=(-1,0)
\) to \(
z_T=(1,0).
\)
As shown in Figure~\ref{fig:sisyphus_mptp}, the simulated transition
trajectories fluctuate around the numerically computed MPTP, showing good
agreement in both the time histories and the phase-plane dynamics.

\subsection{Stochastic damped mechanical system}
The two examples above are two-dimensional degenerate systems. We now consider
a three-dimensional degenerate system driven by multiplicative noise. This
example is inspired by the stochastic damped mechanical system studied in
\cite{zhan_numerical_2023}, which belongs to the class of stochastic contact
systems. In contrast to the Stratonovich model driven by a common noise source
in \cite{zhan_numerical_2023}, we consider here a modified It\^{o} formulation
with two independent standard Brownian motions:
\begin{equation}
\left\{
\begin{aligned}
    \mathrm{d}q_t
    &=
    p_t\,\mathrm{d}t,\\
    \mathrm{d}p_t
    &=
    -\left[
        V'(q_t)+\alpha s_t p_t
    \right]\mathrm{d}t
    -\cos q_t\,\mathrm{d}W_{1,t},\\
    \mathrm{d}s_t
    &=
    \left[
        \frac{1}{2}p_t^2
        -V(q_t)
        -\frac{1}{2}\alpha s_t^2
    \right]\mathrm{d}t
    -\sin q_t\,\mathrm{d}W_{2,t}.
\end{aligned}
\right.
\label{eq:modified_stochastic_contact_system}
\end{equation}
Here \(W_{1,t}\) and \(W_{2,t}\) are independent standard Brownian motions, and
\(\alpha>0\) is the damping parameter.

To put \eqref{eq:modified_stochastic_contact_system} into the framework of
Theorem~\ref{theorem3.2}, we set
\(
    X=q,
    Y=(p,s)^T.
\)
Then the system can be written as
\[
\left\{
\begin{aligned}
    \mathrm{d}X_t
    &=
    g(X_t,Y_t)\,\mathrm{d}t,\\
    \mathrm{d}Y_t
    &=
    f(X_t,Y_t)\,\mathrm{d}t
    +
    \sigma(X_t,Y_t)\,\mathrm{d}W_t,
\end{aligned}
\right.
\]
where \(W_t=(W_{1,t},W_{2,t})^T\),
\(
    g(q,p,s)=p,
\)
\(
    f(q,p,s)
    =
    \begin{pmatrix}
        -V'(q)-\alpha sp\\[1mm]
        \dfrac{1}{2}p^2
        -V(q)
        -\dfrac{1}{2}\alpha s^2
    \end{pmatrix},
\)
and
\(
    \sigma(q,p,s)
    =
    \sigma(q)
    =
    \begin{pmatrix}
        -\cos q & 0\\
        0 & -\sin q
    \end{pmatrix}.
\)
The system is degenerate in the full \((q,p,s)\)-space, since the noise does
not act directly on the \(q\)-component.

We work on a local domain
\(
    \mathcal D=I\times\mathbb R^2,
\)
where \(I\subset\mathbb R\) is an open interval satisfying
\(
    \operatorname{dist}
    \left(
        \overline I,\frac{\pi}{2}\mathbb Z
    \right)>0.
\)
Since \(\frac{\pi}{2}\mathbb Z\) contains all zeros of both \(\sin q\) and
\(\cos q\), this condition ensures that the admissible range of \(q\) stays a
positive distance away from the singular set of \(\sigma(q)^{-1}\). In
particular, there exists a constant \(c_I>0\) such that
\(
    |\sin q|\ge c_I,
    |\cos q|\ge c_I,
q\in I.
\)
No additional restriction is imposed on the \((p,s)\)-variables. Hence
\(\sigma(q)\) is smooth and uniformly invertible on \(\mathcal D\).
The corresponding partial Lamperti transformation is
\[
    \Gamma(q,p,s)
    =
    \begin{pmatrix}
        -p\sec q\\[1mm]
        -s\csc q
    \end{pmatrix}.
\]
Indeed,
\(
    \nabla_{(p,s)}\Gamma(q,p,s)
    =
    \begin{pmatrix}
        -\sec q & 0\\
        0 & -\csc q
    \end{pmatrix}
    =
    \sigma(q)^{-1}.
\)
For each fixed \(q\in I\), the map
\((p,s)\mapsto \Gamma(q,p,s)\) is a \(C^\infty\)-diffeomorphism from
\(\mathbb R^2\) onto \(\mathbb R^2\), with inverse
\(
    (p,s)^T
    =
    v(q,\tilde p,\tilde s)
    =
    \begin{pmatrix}
        -\tilde p\cos q\\[1mm]
        -\tilde s\sin q
    \end{pmatrix}.
\)
Moreover, since \(\sin q\) and \(\cos q\) are bounded away from zero on \(I\),
this family of diffeomorphisms is smooth and uniformly non-degenerate on
\(\mathcal D\). Therefore, the assumptions of Theorem~\ref{theorem3.2} are
satisfied on \(\mathcal D\).

According to Theorem~\ref{theorem3.2}, for the quadratic potential
\(
    V(q)=\frac{1}{2}q^2,
\)
the OM function of system~\eqref{eq:modified_stochastic_contact_system} is
given by
\begin{align*}
    \operatorname{OM}(\phi,\dot{\phi})
    =
    &
    \frac{1}{2}
    \sec^2\phi_1
    \left(
        \dot{\phi}_2+\phi_1+\alpha\phi_2\phi_3
    \right)^2
   +
    \frac{1}{2}
    \csc^2\phi_1
    \left(
        \dot{\phi}_3
        -
        \frac{1}{2}\phi_2^2
        +
        \frac{1}{2}\phi_1^2
        +
        \frac{1}{2}\alpha\phi_3^2
    \right)^2
    \\
    &+
    \phi_2\tan\phi_1
    -
    \frac{1}{2}\phi_2\cot\phi_1
    -
    \alpha\phi_3 .
\end{align*}
Here \(\phi=(\phi_1,\phi_2,\phi_3)^T\) corresponds to \((q,p,s)^T\), and the
formula is applied to admissible paths satisfying
\(
    \phi(t)\in\mathcal D,
    t\in[0,T].
\)

For notational convenience, set
\(
    M(q,p,s)
    =
    2p\tan q-p\cot q-2\alpha s .
\)
By introducing the control
\(
    \theta(t)=(\theta_1(t),\theta_2(t))^T
\)
through
\(
    \dot{Y}(t)
    =
    f(\phi(t))
    +
    \sigma(\phi_1(t))\theta(t),
\)
the corresponding MPTP problem can be reformulated as the following optimal
control problem:
\begin{equation}
\left\{
\begin{aligned}
    \min_{\theta,\phi}\quad
    &
    \int_0^T
    \left[
        \frac{1}{2}|\theta_1(t)|^2
        +
        \frac{1}{2}|\theta_2(t)|^2
        +
        \frac{1}{2}M(\phi_1(t),\phi_2(t),\phi_3(t))
    \right]\mathrm{d}t,
    \\
    \text{subject to}\quad
    &
    \dot{\phi}_1(t)=\phi_2(t),
    \\
    &
    \dot{\phi}_2(t)
    =
    -\phi_1(t)-\alpha \phi_2(t)\phi_3(t)
    -\cos\phi_1(t)\,\theta_1(t),
    \\
    &
    \dot{\phi}_3(t)
    =
    \frac{1}{2}\phi_2^2(t)
    -
    \frac{1}{2}\phi_1^2(t)
    -
    \frac{1}{2}\alpha\phi_3^2(t)
    -
    \sin\phi_1(t)\,\theta_2(t),
    \\
    &
    \phi(0)=z_0,
    \qquad
    \phi(T)=z_T,
    \qquad
    \phi(t)\in\mathcal D,\quad t\in[0,T].
\end{aligned}
\right.
\label{eq:contact_optimal_control}
\end{equation}
Here \(z_0,z_T\in\mathcal D\) denote the prescribed initial and terminal states.

In this example, the diffusion matrix depends only on the non-noisy variable
\(q\), and therefore no additional shift term appears in the quadratic control cost. In the numerical experiment, we choose
\(
    I=(0.3,1.2),
    \mathcal D=I\times\mathbb R^2,
\)
and compute the transition path between
\(
    z_0=(0.5,0,0)^T,
    z_T=(1.0,0,0)^T.
\)
The computed path satisfies \(\phi_1(t)\in I\) for all \(t\in[0,T]\), and
therefore remains in the admissible local domain where the diffusion matrix is
smooth and uniformly invertible.

\begin{figure}[H]
    \centering
    \includegraphics[width=0.8\linewidth]{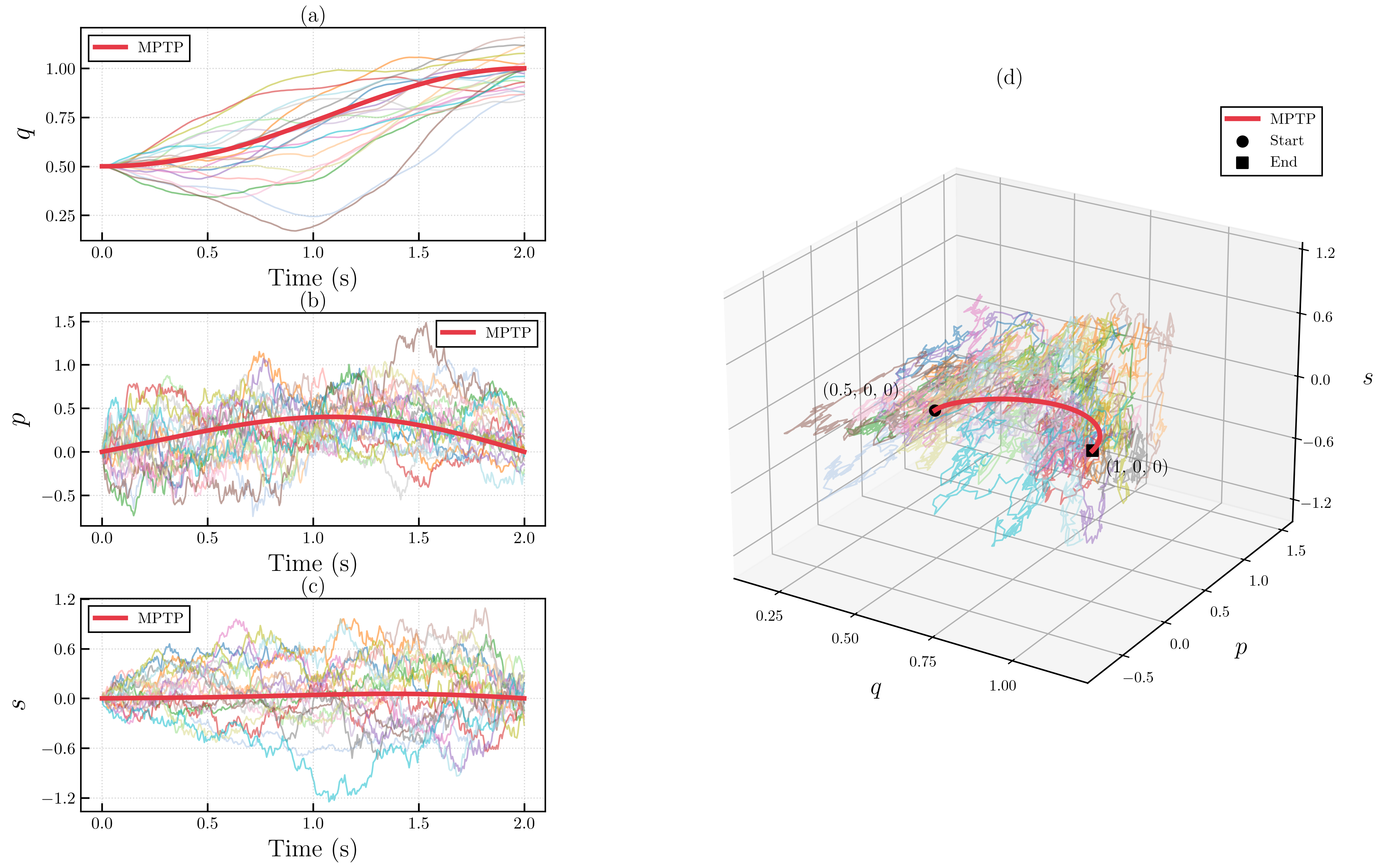}
    \caption{
    Numerical MPTP and simulated trajectories for the stochastic damped
    mechanical system.
    }
    \label{fig:Mfi}
\end{figure}

For \(\alpha=0.1\), Figure~\ref{fig:Mfi} compares the MPTP computed by the
proposed method with sample trajectories simulated by the Euler--Maruyama
method. The numerical results indicate that the proposed method remains
effective for computing MPTPs in stochastic dynamical systems driven by
degenerate multiplicative noise. They also demonstrate that the
optimal-control-based neural network approach provides a practical tool for
computing MPTPs in more complex stochastic systems.

\section{Conclusion}\label{sec5}

In this paper, we investigated the most probable transition paths (MPTPs) for
noise-induced transitions in degenerate stochastic systems driven by
multiplicative Gaussian noise. For a class of degenerate stochastic differential
equations with state-dependent diffusion, we derived an explicit
Onsager--Machlup (OM) action functional. The derivation is based on a
Lamperti-type transformation applied to the noisy component, combined with the
Girsanov theorem, which reduces the problem to an additive-noise formulation in
the transformed variables. The resulting functional extends existing OM
formulations for nondegenerate diffusions and degenerate additive-noise systems,
and provides a variational characterization of MPTPs under degenerate
multiplicative noise.

To compute transition pathways, we reformulated the minimization of the OM
action functional as an optimal control problem, thereby avoiding the direct
solution of the associated Hamilton--Pontryagin boundary value problem. A
fully connected neural-network method was then developed to approximate the
optimal control and the corresponding transition path. Numerical experiments on
a kinetic Langevin system, a Sisyphus-cooling system, and a
stochastic damped mechanical system show good agreement with available
analytical solutions and with transition behaviors observed in stochastic
simulations. These results demonstrate the effectiveness of the proposed method in capturing dominant transition routes under degenerate multiplicative noise.

\section*{Appendix. Proof of equation~\eqref{conclusion1}}

\setcounter{equation}{0}
\renewcommand{\theequation}{A.\arabic{equation}}

The proof consists of two steps. We first
establish the supremum-norm equivalence between \(\tilde Y-\phi_2\) and the
Brownian motion \(W\). Then, using the equation for the \(\tilde X\)-component
and the global Lipschitz continuity of \(g\), we control
\(\tilde X-\phi_1\) by \(\tilde Y-\phi_2\) in the supremum norm. This yields the
corresponding supremum-norm equivalence between \(\tilde Z-\phi\) and \(W\).

Let
\[
    \overline Y_t=\tilde Y_t-\phi_2(t),\qquad
    \overline X_t=\tilde X_t-\phi_1(t),\qquad
    \overline Z_t=(\overline X_t,\overline Y_t),
\]
where \(\tilde Z=(\tilde X,\tilde Y)\) is given by \eqref{ASDE}. For a continuous
process \(U\) on \([0,T]\), we write
\(
    \|U\|:=\sup_{0\leq t\leq T}|U_t|.
\)
The generic constants below may depend on \(T\), on the Lipschitz constant of
\(g\), and on the uniform bounds for \(\sigma\), \(\sigma^{-1}\), and
\(\dot\sigma\), but are independent of the Brownian path \(W\).

\medskip
\noindent\textbf{Step 1: Supremum-norm equivalence between \(\overline Y\) and \(W\).}

Since \(\sigma(t)\) is uniformly nondegenerate and continuous on \([0,T]\),
there exist constants \(0<\lambda_*\leq \lambda^*<\infty\) such that
\[
    \lambda_* |x|
    \leq
    |\sigma(t)x|
    \leq
    \lambda^* |x|,
    \qquad
    x\in\mathbb R^m,\quad t\in[0,T].
\]
Moreover, since \(\dot\sigma\) is bounded on \([0,T]\), there exists
\(L_\sigma>0\) such that
\[
    \|\dot\sigma(t)\|_{\mathrm{op}}\leq L_\sigma,
    \qquad t\in[0,T].
\]
Here \(\|\cdot\|_{\mathrm{op}}\) denotes the induced operator norm for matrices, namely
\(
\|A\|_{\mathrm{op}}=\sup_{|x|=1}|Ax|.
\) From the auxiliary equation \eqref{ASDE}, we have
\[
    \overline Y_t=\int_0^t \sigma(s)\,\mathrm dW_s .
\]
Since \(\sigma\) is deterministic and continuously differentiable, the
integration-by-parts formula gives
\[
    \overline Y_t
    =
    \sigma(t)W_t
    -
    \int_0^t \dot\sigma(s)W_s\,\mathrm ds .
\]
Hence
\[
    |\overline Y_t|
    \leq
    \lambda^* |W_t|
    +
    L_\sigma\int_0^t |W_s|\,\mathrm ds
    \leq
    \left(\lambda^*+TL_\sigma\right)\|W\|.
\]
Taking the supremum over \(t\in[0,T]\), we obtain
\[
    \|\overline Y\|
    \leq
    \left(\lambda^*+TL_\sigma\right)\|W\|.
\]
Conversely, from
\[
    \sigma(t)W_t
    =
    \overline Y_t
    +
    \int_0^t \dot\sigma(s)W_s\,\mathrm ds ,
\]
we get
\[
    \lambda_* |W_t|
    \leq
    |\overline Y_t|
    +
    L_\sigma\int_0^t |W_s|\,\mathrm ds .
\]
Therefore,
\[
    |W_t|
    \leq
    \frac{1}{\lambda_*}\|\overline Y\|
    +
    \frac{L_\sigma}{\lambda_*}
    \int_0^t |W_s|\,\mathrm ds .
\]
By Gronwall's inequality,
\[
    |W_t|
    \leq
    \frac{1}{\lambda_*}
    \exp\left(\frac{L_\sigma t}{\lambda_*}\right)
    \|\overline Y\|,
    \qquad t\in[0,T].
\]
Taking the supremum over \(t\in[0,T]\), we obtain
\[
    \|W\|
    \leq
    \frac{1}{\lambda_*}
    \exp\left(\frac{L_\sigma T}{\lambda_*}\right)
    \|\overline Y\|.
\]
Consequently, there exist constants \(a_1,a_2>0\), independent of \(W\), such
that
\begin{equation}
\label{eq:Y-W-sup}
    a_1\|W\|
    \leq
    \|\overline Y\|
    \leq
    a_2\|W\| .
\end{equation}

\medskip
\noindent\textbf{Step 2: Supremum-norm equivalence between \(\overline Z\) and \(W\).}

Since \(\overline Z_t=(\overline X_t,\overline Y_t)\), we immediately have
\[
    \|\overline Y\|\leq \|\overline Z\|.
\]
It remains to control \(\overline X\) by \(\overline Y\). By the equation for
\(\tilde X\) in \eqref{ASDE} and the constraint satisfied by the reference path,
namely
\[
    \phi_1(t)=x_0+\int_0^t g(s,\phi_1(s),\phi_2(s))\,\mathrm ds,
\]
we obtain
\[
    \overline X_t
    =
    \int_0^t
    \left[
        g(s,\tilde X_s,\tilde Y_s)
        -
        g(s,\phi_1(s),\phi_2(s))
    \right]
    \mathrm ds .
\]
Since \(g\in C_b^{0,1}\), it is globally Lipschitz in \((x,y)\), uniformly in
\(t\). Thus, for some constant \(L_g>0\),
\[
    |\overline X_t|
    \leq
    L_g\int_0^t
    \left(
        |\overline X_s|
        +
        |\overline Y_s|
    \right)
    \mathrm ds .
\]
Hence
\(
    |\overline X_t|
    \leq
    L_g\int_0^t |\overline X_s|\,\mathrm ds
    +
    L_gT\|\overline Y\|.
\)
By Gronwall's inequality,
\[
    |\overline X_t|
    \leq
    L_gT e^{L_gT}\|\overline Y\|,
    \qquad t\in[0,T].
\]
Taking the supremum over \(t\in[0,T]\), we find that there exists a constant
\(K_2>0\), independent of \(W\), such that
\begin{equation}
\label{eq:X-controlled-by-Y}
    \|\overline X\|
    \leq
    K_2\|\overline Y\|.
\end{equation}
Consequently,
\[
    \|\overline Z\|
    \leq
    \|\overline X\|+\|\overline Y\|
    \leq
    (K_2+1)\|\overline Y\|.
\]
Together with \(\|\overline Y\|\leq \|\overline Z\|\), this gives
\begin{equation}
\label{eq:Z-Y-sup}
    \|\overline Y\|
    \leq
    \|\overline Z\|
    \leq
    (K_2+1)\|\overline Y\|.
\end{equation}
Combining \eqref{eq:Z-Y-sup} with \eqref{eq:Y-W-sup}, we conclude that there
exist constants \(c_1,c_2>0\), independent of \(W\), such that
\begin{equation}
\label{eq:Z-W-sup}
    c_1\|W\|
    \leq
    \|\overline Z\|
    \leq
    c_2\|W\|.
\end{equation}
Since \(\overline Z=\tilde Z-\phi\), equation~\eqref{conclusion1} follows.
\hfill \(\qedsymbol\)


\section*{Acknowledgements}
The research was partially supported by the Natural Science Basic Research
Program of Shaanxi Province (Program No.~2026JC-YBQN-0096), the Shaanxi
Fundamental Science Research Project for Mathematics and Physics
(Grant No.~25JSQ042), the Fundamental Research Funds for the Central
Universities (Grant No.~xzy012025071), the National Natural Science
Foundation of China (Grant Nos.~12101484 and 12141107), the Guangdong
Provincial Key Laboratory of Mathematical and Neural Dynamical Systems
(No.~2024B1212010004), the Cross-Disciplinary Research Team on Data Science
and Intelligent Medicine (No.~2023KCXTD054), and the Guangdong--Dongguan
Joint Research Grant (No.~2023A1515140016).

\bibliographystyle{alpha}
\bibliography{Refs}

@book{Ikeda2014stochastic,
  title={Stochastic differential equations and diffusion processes},
  author={Ikeda, N. and Watanabe, S.},
  volume={24},
  year={2014},
  publisher={Elsevier}
}

@article{cheng_most_2019,
  title={Most probable transition pathways and maximal likely trajectories in a genetic regulatory system},
  author={Cheng, X. and Wang, H. and Wang, X. and Duan, J. and Li, X.},
  journal={Phys. A},
  volume={531},
  pages={121779},
  year={2019},
  publisher={Elsevier}
}

@article{zheng_transitions_2016,
  title={Transitions in a genetic transcriptional regulatory system under {L}{\'e}vy motion},
  author={Zheng, Y. and Serdukova, L. and Duan, J. and Kurths, J.},
  journal={Sci. Rep.},
  volume={6},
  number={1},
  pages={29274},
  year={2016},
  publisher={Nature Publishing Group UK London}
}

@article{liu_fluctuations-induced_2004,
  title={Fluctuations-induced switch in the gene transcriptional regulatory system},
  author={Liu, Q. and Jia, Y.},
  journal={Phys. Rev. E},
  volume={70},
  number={4},
  pages={041907},
  year={2004},
  publisher={APS}
}

@article{wei_most_2021,
  title={Most probable transitions from metastable to oscillatory regimes in a carbon cycle system},
  author={Wei, W. and Hu, J. and Chen, J. and Duan, J.},
  journal={Chaos},
  volume={31},
  number={12},
  pages={121102},
  year={2021},
  publisher={AIP Publishing}
}

@article{onsager_fluctuations_1953,
  title={Fluctuations and irreversible processes},
  author={Onsager, L. and Machlup, S.},
  journal={Phys. Rev.},
  volume={91},
  number={6},
  pages={1505--1512},
  year={1953},
  publisher={APS}
}

@article{tisza_fluctuations_1957,
  title={Fluctuations and irreversible thermodynamics},
  author={Tisza, L. and Manning, I.},
  journal={Phys. Rev.},
  volume={105},
  number={6},
  pages={1695--1705},
  year={1957},
  publisher={APS}
}

@article{stratonovich1971probability,
  title={On the probability functional of diffusion processes},
  author={Stratonovich, R. L.},
  journal={Selected Trans. in Math. Stat. Prob},
  volume={10},
  pages={273--286},
  year={1971}
}

@article {capitaine_onsager-machlup_1995,
    AUTHOR = {Capitaine, M.},
     TITLE = {Onsager--{M}achlup functional for some smooth norms on {W}iener
              space},
   JOURNAL = {Probab. Theory Related Fields},
    VOLUME = {102},
      YEAR = {1995},
    NUMBER = {2},
     PAGES = {189--201},
      ISSN = {0178-8051,1432-2064},
   MRCLASS = {60J60 (60G15 60H10 60J65)},
  MRNUMBER = {1337251},
MRREVIEWER = {Ofer\ Zeitouni},
       DOI = {10.1007/BF01213388},
       URL = {https://doi.org/10.1007/BF01213388},
}

@article {durr_onsager-machlup_1978,
    AUTHOR = {D\"urr, D. and Bach, A.},
     TITLE = {The {O}nsager--{M}achlup function as {L}agrangian for the most
              probable path of a diffusion process},
   JOURNAL = {Comm. Math. Phys.},
    VOLUME = {60},
      YEAR = {1978},
    NUMBER = {2},
     PAGES = {153--170},
      ISSN = {0010-3616,1432-0916},
   MRCLASS = {82.60},
  MRNUMBER = {489608},
MRREVIEWER = {Burton\ H.\ Voorhees},
       URL = {http://projecteuclid.org/euclid.cmp/1103904077},
}

@article {chao_onsagermachlup_2019,
    AUTHOR = {Chao, Y. and Duan, J.},
     TITLE = {The {O}nsager--{M}achlup function as {L}agrangian for the most
              probable path of a jump-diffusion process},
   JOURNAL = {Nonlinearity},
    VOLUME = {32},
      YEAR = {2019},
    NUMBER = {10},
     PAGES = {3715--3741},
      ISSN = {0951-7715,1361-6544},
   MRCLASS = {60G17 (60J75 82C40)},
  MRNUMBER = {4002397},
       DOI = {10.1088/1361-6544/ab248b},
       URL = {https://doi.org/10.1088/1361-6544/ab248b},
}

@article{horsthemke_onsager-machlup_1975,
  title={Onsager--{M}achlup function for one dimensional nonlinear diffusion processes},
  author={Horsthemke, W and Bach, A},
  journal={Z. Phys. B},
  volume={22},
  number={2},
  pages={189--192},
  year={1975},
  publisher={Springer}
}

@article {bach_functionals_1977,
    AUTHOR = {Bach, A. and D\"urr, D. and Stawicki, B.},
     TITLE = {Functionals of paths of a diffusion process and the
              {O}nsager--{M}achlup function},
   JOURNAL = {Z. Phys. B},
    VOLUME = {26},
      YEAR = {1977},
    NUMBER = {2},
     PAGES = {191--193},
      ISSN = {0722-3277},
   MRCLASS = {82.60},
  MRNUMBER = {432125},
MRREVIEWER = {J.\ Fritz},
       DOI = {10.1007/BF01325272},
       URL = {https://doi.org/10.1007/BF01325272},
}

@article{hu_transition_2021,
  title={Transition pathways for a class of high dimensional stochastic dynamical systems with {L}{\'e}vy noise},
  author={Hu, J. and Chen, J.},
  journal={Chaos},
  volume={31},
  number={6},
  pages={063138},
  year={2021},
  publisher={AIP Publishing}
}

@article {WOS:001473901800007,
    AUTHOR = {Huang, Y. and Zhou, X. and Duan, J.},
     TITLE = {Probability flow approach to the {O}nsager--{M}achlup
              functional for jump-diffusion processes},
   JOURNAL = {SIAM J. Appl. Math.},
    VOLUME = {85},
      YEAR = {2025},
    NUMBER = {2},
     PAGES = {524--547},
      ISSN = {0036-1399,1095-712X},
   MRCLASS = {60H10 (37H05 60J60 82C31)},
  MRNUMBER = {4876543},
       DOI = {10.1137/24M1689259},
       URL = {https://doi.org/10.1137/24M1689259},
}

@article{WOS:001514271800001,
  title={Transition pathways for a class of degenerate stochastic dynamical systems with {L}{\'e}vy noise},
  author={Chao, Y. and Wei, P.},
  journal={Commun. Nonlinear Sci. Numer. Simul.},
  volume={150},
  number={12},
  pages={109034},
  year={2025},
  publisher={Elsevier}
}

@article{aihara_mortensen_1999,
  title={On the {M}ortensen equation for maximum likelihood state estimation},
  author={Aihara, S. I. and Bagchi, A.},
  journal={IEEE Trans. Automat. Control},
  volume={44},
  number={10},
  pages={1955--1961},
  year={1999},
  publisher={IEEE}
}

@article{liu_onsagermachlup_2024,
  title={The Onsager--Machlup action functional for degenerate stochastic differential equations in a class of norms},
  author={Liu, S. and Gao, H.},
  journal={Statist. Probab. Lett.},
  volume={206},
  pages={110009},
  year={2024},
  publisher={Elsevier}
}

@article{li2021machine,
  title={Machine learning framework for computing the most probable paths of stochastic dynamical systems},
  author={Li, Y. and Duan, J. and Liu, X.},
  journal={Phys. Rev. E},
  volume={103},
  number={1},
  pages={012124},
  year={2021},
  publisher={APS}
}

@article{soons_optimal_2024,
	title = {Optimal Transition Paths for {AMOC} Collapse and Recovery in a Stochastic Box Model},
	volume = {54},
	copyright = {http://www.ametsoc.org/PUBSReuseLicenses},
	issn = {0022-3670, 1520-0485},
	url = {https://journals.ametsoc.org/view/journals/phoc/54/12/JPO-D-23-0234.1.xml},
	doi = {10.1175/JPO-D-23-0234.1},
	language = {en},
	number = {12},
	urldate = {2026-04-29},
	journal = {J. Phys. Oceanogr.},
	author = {Soons, J. and Grafke, T. and Dijkstra, H. A.},
	year = {2024},
	pages = {2537--2552},
}

@article {doi:10.1137/19M1288802,
    AUTHOR = {Nakamura-Zimmerer, T. and Gong, Q. and Kang, W.},
     TITLE = {Adaptive deep learning for high-dimensional
              {H}amilton-{J}acobi-{B}ellman equations},
   JOURNAL = {SIAM J. Sci. Comput.},
    VOLUME = {43},
      YEAR = {2021},
    NUMBER = {2},
     PAGES = {A1221--A1247},
      ISSN = {1064-8275,1095-7197},
   MRCLASS = {49N35 (49K15 49L20 68T07 93C15)},
  MRNUMBER = {4241190},
       DOI = {10.1137/19M1288802},
       URL = {https://doi.org/10.1137/19M1288802},
}

@article {chen_detecting_2024,
    AUTHOR = {Chen, J. and Gao, T. and Li, Y. and Duan, J.},
     TITLE = {Detecting the most probable transition pathway based on
              optimal control theory},
   JOURNAL = {Appl. Math. Model.},
    VOLUME = {127},
      YEAR = {2024},
     PAGES = {217--236},
      ISSN = {0307-904X,1872-8480},
   MRCLASS = {49N90 (49K15 60H10 65K10)},
  MRNUMBER = {4678649},
       DOI = {10.1016/j.apm.2023.12.008},
       URL = {https://doi.org/10.1016/j.apm.2023.12.008},
}

@article{borell1977gaussian,
  title   = {Gaussian measures and the isoperimetric inequality},
  author  = {Borell, C.},
  journal = {Ann. Inst. H. Poincaré Sect. B},
  volume  = {13},
  number  = {3},
  pages   = {245--258},
  year    = {1977}
}

@incollection {ledoux1996isoperimetry,
    AUTHOR = {Ledoux, M.},
     TITLE = {Isoperimetry and {G}aussian analysis},
 BOOKTITLE = {Lectures on probability theory and statistics
              ({S}aint-{F}lour, 1994)},
    SERIES = {Lecture Notes in Math.},
    VOLUME = {1648},
     PAGES = {165--294},
 PUBLISHER = {Springer, Berlin},
      YEAR = {1996},
      ISBN = {3-540-62055-9},
   MRCLASS = {60-02 (60E15 60F10 60G15)},
  MRNUMBER = {1600888},
MRREVIEWER = {Qi\ Man\ Shao},
       DOI = {10.1007/BFb0095676},
       URL = {https://doi.org/10.1007/BFb0095676},
}

@article{sheppzeitouni1993,
  title   = {A note on the conditional exponential moments of {W}iener integrals},
  author  = {Shepp, L. and Zeitouni, O.},
  journal = {Ann. Probab.},
  volume  = {21},
  number  = {2},
  pages   = {1124--1134},
  year    = {1993}
}

@article {lan_new_2014,
    AUTHOR = {Lan, G. and Wu, J.-L.},
     TITLE = {New sufficient conditions of existence, moment estimations and
              non confluence for {SDE}s with non-{L}ipschitzian
              coefficients},
   JOURNAL = {Stochastic Process. Appl.},
    VOLUME = {124},
      YEAR = {2014},
    NUMBER = {12},
     PAGES = {4030--4049},
      ISSN = {0304-4149,1879-209X},
   MRCLASS = {60H10},
  MRNUMBER = {3264438},
MRREVIEWER = {Paul\ Raynaud de Fitte},
       DOI = {10.1016/j.spa.2014.07.010},
       URL = {https://doi.org/10.1016/j.spa.2014.07.010},
}

@article {wei_optimal_2022,
    AUTHOR = {Wei, W. and Gao, T. and Chen, X. and Duan, J.},
     TITLE = {An optimal control method to compute the most likely
              transition path for stochastic dynamical systems with jumps},
   JOURNAL = {Chaos},
    VOLUME = {32},
      YEAR = {2022},
    NUMBER = {5},
     PAGES = {051102},
      ISSN = {1054-1500,1089-7682},
   MRCLASS = {60H30 (34F05 60F10)},
  MRNUMBER = {4422807},
       DOI = {10.1063/5.0093924},
       URL = {https://doi.org/10.1063/5.0093924},
}

@inproceedings{baxendale_lyapunov_2003,
  title={Lyapunov exponents and stability for the stochastic {D}uffing-van der {P}ol oscillator},
  author={Baxendale, P. H.},
  booktitle={IUTAM Symposium on Nonlinear Stochastic Dynamics: Proceedings of the IUTAM Symposium},
  pages={125--135},
  year={2003},
  organization={Springer}
}

@article{zhan_numerical_2023,
  title={Numerical integration of stochastic contact {H}amiltonian systems via stochastic {H}erglotz variational principle},
  author={Zhan, Q. and Duan, J. and Li, X. and Li, Y.},
  journal={Phys. Scripta},
  volume={98},
  number={5},
  pages={055211},
  year={2023},
  publisher={IOP Publishing}
}

@techreport{moller2010state,
  author      = {M{\o}ller, J. K. and Madsen, H.},
  title       = {From State Dependent Diffusion to Constant Diffusion in Stochastic Differential Equations by the  {L}amperti Transform},
  institution = {Technical University of Denmark, DTU Informatics},
  number= {IMM-Technical Report-2010-16},
  year= {2010}
}

@article{Falasco2022,
  title={Generalized virial equation for nonlinear multiplicative Langevin dynamics: application to laser-cooled atoms},
  author={Falasco, G. and Barkai, E. and Baiesi, M.},
  journal={Phys. Rev. E},
  volume={105},
  number={2},
  pages={024143},
  year={2022},
  publisher={APS}
}
\end{document}